\documentclass[review,onefignum,onetabnum]{siamart171218}

\usepackage{lipsum}
\usepackage{amsfonts}
\usepackage{graphicx}
\usepackage{epstopdf}
\ifpdf
  \DeclareGraphicsExtensions{.eps,.pdf,.png,.jpg}
\else
  \DeclareGraphicsExtensions{.eps}
\fi

\newsiamremark{remark}{Remark}
\newsiamremark{hypothesis}{Hypothesis}
\crefname{hypothesis}{Hypothesis}{Hypotheses}
\newsiamthm{claim}{Claim}

\title{
Stochastic Three Points Method \\ for Unconstrained Smooth Minimization}

 \author{
El Houcine Bergou\thanks{King  Abdullah  University  of  Science  and  Technology  (KAUST),  Thuwal,  Saudi  Arabia. MaIAGE, INRA, Universit\'e Paris-Saclay, 78350 Jouy-en-Josas, France
 ({\tt elhoucine.bergou@inra.fr}). This author received support from the AgreenSkills+ fellowship programme which has received funding from the EU's Seventh Framework Programme under grant agreement No FP7-609398 (AgreenSkills+ contract).
}
 \and Eduard Gorbunov\thanks{Moscow Institute of Physics and Technology (MIPT), Moscow, Russian Federation
({\tt eduard.gorbunov@phystech.edu}).
}
  \and Peter Richt\'{a}rik\thanks{
  King  Abdullah  University  of  Science  and  Technology  (KAUST),  Thuwal,  Saudi  Arabia. University  of  Edinburgh, Edinburgh, United Kingdom. Moscow Institute of Physics and Technology (MIPT), Moscow, Russian Federation ({\tt peter.richtarik@kaust.edu.sa}).}
}

\usepackage{amsopn}

\usepackage{tikz}
\usepackage{subfigure} 
\usepackage{amssymb}
\usepackage{pifont}

\newcommand{\circledOne}{\text{\ding{172}}}
\newcommand{\circledTwo}{\text{\ding{173}}}

\newcommand{\cD}{\cal D}
 
\newcommand{\e}{\varepsilon}
\newcommand{\Exp}{\mathbf{E}}

\newcommand{\R}{\mathbb{R}}
\newcommand{\eqdef}{\stackrel{\text{def}}{=}}
\newcommand{\ve}[2]{\left\langle #1 , #2 \right\rangle}
\def\<#1,#2>{\left\langle #1,#2\right\rangle}

\newtheorem{ass}[theorem]{Assumption}
\newtheorem{lem}[theorem]{Lemma}
\newtheorem{thm}[theorem]{Theorem}

\ifpdf
\hypersetup{
  pdftitle={Stochastic Three Points Method for Unconstrained Smooth Minimization},
  pdfauthor={}
}
\fi

\begin{document}

\maketitle

\begin{abstract}
  In this paper we consider the unconstrained minimization  problem of a smooth function in $\R^n$ in a setting where only function evaluations are possible. We design a novel randomized derivative-free algorithm --- the {\em stochastic three  points (STP)} method --- and analyze its iteration complexity. At each iteration, \texttt{STP} generates a random search direction according to a certain fixed probability law. Our assumptions on this law are very mild: roughly speaking,  all laws which do not concentrate all measure on any halfspace passing through the origin will work. For instance, we allow for the uniform distribution on the sphere and also  distributions that concentrate all measure on a positive spanning set. 
 
 Although our approach is designed to not use explicitly derivatives, it covers some first order methods. For instance if the probability law is chosen to be the Dirac distribution concentrated at the sign of the gradient then \texttt{STP} recovers the Signed Gradient Descent method. If the probability law is the uniform distribution on the coordinates of the gradient then \texttt{STP} recovers the Coordinate Descent Method.
 
 Given a current iterate $x$, \texttt{STP} compares  the objective function  at three  points: $x$, $x+\alpha s$ and $x-\alpha s$, where $\alpha>0$ is a stepsize parameter and $s$ is the random search direction. The best of these three points is the next iterate.  We analyze the method \texttt{STP} under several stepsize selection schemes (fixed, decreasing, estimated through finite differences, etc). 
 
The complexity of \texttt{STP} depends on the probability law via a simple characteristic closely related to the cosine measure which is used in the analysis of deterministic direct search (\texttt{DDS}) methods. Unlike in \texttt{DDS}, where  $O(n)$ ($n$ is the dimension of $x$) function evaluations must be performed in each iteration in the worst case, our method only requires two new function evaluations per iteration. Consequently, while \texttt{DDS} depends quadratically on $n$, our method depends linearly on $n$. In particular, in the nonconvex case, \texttt{STP} needs  $O(n\e^{-2})$ function evaluations to find a point at which the gradient of the objective function is below $\e$, in expectation. In the convex case, the complexity is $O(n \e^{-1})$. In the strongly convex case \texttt{STP} converges linearly, meaning that the complexity is $O\left(n \log\left({\e}^{-1}\right)\right)$. 
\end{abstract}

%

\section{Introduction}

In this paper we consider the problem 
\begin{equation}\label{eq:mainP}
\min_{x\in \R^n} f(x),
\end{equation}
where $f: \R^n \rightarrow \R$ is a given  smooth objective function. We assume that we do not have access to the derivatives of $f$ and  only have access to a function evaluation oracle. In other words, we assume that we work in the Derivative-Free Optimization (DFO) setting  \cite{Conn_2009}. 
Optimization problems of this type appear in many industrial applications where usually the objective function is
evaluated through a computer simulation process, and therefore derivatives cannot be directly evaluated; e.g., shape
optimization in fluid-dynamics problems  \cite{Allaire_2001, Haslinger_2003, Mohammadi_2001}.

 
 Direct search methods of directional type \cite{Kolda_2003, Conn_2009} are a popular class of methods for DFO  and are  among the first algorithms proposed in numerical optimization \cite{Matyas_1965}. 
These methods are characterized by  evaluating the objective function  over a number of (typically predetermined and fixed) directions to ensure descent using a sufficiently small stepsize. The directions are typically required to form a {\em positive spanning set} (i.e.\ a set of vectors whose conic hull is  $\R^n$) in order to make sure that each point in $\R^n$ (and hence also the optimal solution) is achievable by a sequence of positive steps from any starting point.


For instance,  the {\em coordinate search} method uses the coordinate (i.e., standard basic) directions, $e_1,e_2,\dots,e_n$, and their negatives, $-e_1,-e_2,\dots,-e_n$ as the set of admissible directions. Clearly, $\{\pm e_i, \;:\; i=1,2,\dots,n\}$ forms a positive spanning set.

\subsection{Stochastic Three Points method}

 In this paper, we study  a very general {\em randomized} variant of direct search methods, which we call {\em Stochastic Three Points} (\texttt{STP}).  

\texttt{STP} depends on two ``parameters'': a distribution / probability law $\cD$ from which we sample directions, and a stepsize selection rule. At  iteration $k$ of \texttt{STP}, we generate a random direction $s_k$ by sampling from ${\cal D}$, and  then choose the next iterate via   \[x_{k+1} = \arg \min \left\{f(x_k + \alpha_k s_k), f(x_k -  \alpha_k s_k),f(x_k) \right\},\] where $\alpha_k>0$ is an appropriately chosen stepsize. That is, we pick $x_{k+1}$ as the best of the three points $x_k+\alpha_k s_k, x_k-\alpha_k s_k$ and $x_k$ in terms of the function values.


We prove for such a scheme, with several different choices of stepsizes, that the number of iterations sufficient
to guarantee that $\min\limits_{k=0,1,\ldots,K}\Exp \left[\|\nabla f(x_k)\|_{\cD}\right] \leq \e$ is  $O(n\e^{-2})$, where $\| \cdot \|_{\cal D}$ is a norm dependent on ${\cal D}$ which we introduce in Section~\ref{sec:SD} and $\Exp \left[\cdot \right]$ is the expectation. This complexity is global since no assumption is made
on the starting point. If the objective function $f$ is convex, then the number of iterations needed 
to get $x_k$ such that $\Exp\left[f(x_k) - f_*\right]\le \varepsilon$ is $O(n \e^{-1})$ where $f_*$ is the optimal value of $f$. If in addition, $f$ is strongly convex, then we have a
global linear rate of convergence. This is an improvement on deterministic direct search (\texttt{DDS}) where the best known complexity bounds depend quadratically on $n$ 
 and the same way as our scheme in $\e$ \cite{KR-DFO2014, Vicente_2013, Vicente_2016}. 
We propose also a parallel version for \texttt{STP}.

Despite our  approach shares similarities with other randomized algorithmic approaches, the differences are significant. In the sixties a random optimization approach was proposed in \cite{Matyas_1965}. It was proposed to sample a point  randomly around the current iterate and move to this new point if it decreases the objective function. This approach was generalized to cover constrained problems in \cite{Baba_1981}. The theoretical and numerical performances of this approach for nonconvex functions was studied in \cite{Dorea_1983, Sarma_1990}. 
 More recently, the works  
in \cite{Diniz_2008} and \cite{Gratton_2015} use  random searching directions, and impose a decrease condition to whether accept the step or reject it, like in \texttt{DDS}. They update the stepsize by increasing it if the step is accepted and decreasing it otherwise.  
 Our approach is different from these frameworks  in the sense that at each iteration we generate a single direction, then we choose the stepsize independently from any decrease condition. 
 In \cite{Gratton_2015}, the authors impose to the search direction some probabilistic property. In fact, they assume that at each iteration their random directions are probabilistic descent conditioned to the past. In other words, at a given iteration, independently from the past with a certain probability at least one of the directions is of descent type. 
 The main result of \cite{Gratton_2015} is the complexity bound $O(r n\e^{-2})$ to drive the gradient norm below $\e$ with high probability, where $r \ge 2$ is the number of the random directions at each iteration. Also \cite{Gratton_2015} do not cover the cases when the objective function is convex or strongly convex. \texttt{STP} method gives similar complexity bound for non-convex problems (with $r=2$).  


More related to our work is the method proposed in \cite{Karmanov_1974a, Karmanov_1974b} for convex problems, where at iteration $k$ the step is updated as follows
\begin{equation*}
\label{eq:xx-method}
x_{k+1} = x_k + \alpha_k u, 
\end{equation*} 
where $u$ is sampled uniformly from the uniform distribution on the unit sphere, and $$\alpha_k = \arg\min_{\alpha \in \R} f(x_k + \alpha u).$$ The latter method was improved  in two ways
by \cite{Stich_2011}. In fact, the proposed method in \cite{Stich_2011} i) allows approximate line search, i.e., $\alpha_k \approx \arg\min_{\alpha \in \R} f(x_k + \alpha u)$, ii) and allows discrete sampling from $\{ \pm e_i, i=1,\ldots,n\}$ instead of sampling from the unit sphere. 
Our approach is different from these methods in the sense that it did not perform any line search approximation to compute the stepsizes, and allows different distributions (which include the uniform distribution over the unit sphere and the discrete sampling from the canonical basis of $\R^n$) to sample the directions.
 The complexity bounds given in these works are worse than those obtained in this paper. 
Another method related to our work is the method discussed in \cite[Section 3.4]{Polyak_1987}, 
a derivative-free approach based on forming an unbiased estimate of the gradient using Gaussian smoothing. The search direction in this method is distributed uniformly over the unit sphere and it is  pre-multiplied by an approximation to the directional derivative along the direction itself. More precisely, this method updates the step at iteration $k$ as follows
\begin{equation}
\label{eq:Polyak-method}
x_{k+1} = x_k - \alpha_k \frac{f(x_k + \mu_k u) - f(x_k)}{ \mu_k} u, 
\end{equation} 
where $\mu_k\in (0,1)$ is the finite differences parameter, $\alpha_k$ is the stepsize,	and $u$ is a random vector distributed uniformly over the unit sphere.
 In this work, 	there is no explicit rules for choosing the parameters and there is 
 no analysis of the worst case complexity. 
 The paper \cite{Nesterov_Spokoiny_2017} proposes other variants of this method
 by changing the way of approximating the directional derivative of $f$ along $u$. Moreover, it gives the worst case complexity analysis of the method (\ref{eq:Polyak-method}). The complexity bounds in \cite{Nesterov_Spokoiny_2017} are similar to those of our \texttt{STP} approach.
 Our approach is different from the method (\ref{eq:Polyak-method}) and its variants porposed in \cite{Nesterov_Spokoiny_2017}, in our approach the search direction  can follows a different distribution from the uniform distribution over the unit sphere. For instance, we allow a distribution that has all its mass concentrated on   a discrete set of vectors -- which makes a direct connection with the (deterministic) direct search methods. Moreover, the proposed stepsizes  in \cite{Nesterov_Spokoiny_2017} depend on the Lipschitz constant of the gradient of the objective function. However, in our approach we proposed some stepsizes which can be easily computed in practice. The extention of the work \cite{Nesterov_Spokoiny_2017} for an uncontrained problem of minimization of a smooth convex function which is only available through noisy observations of its values were studied in the recent work \cite{Dvurechensky_Gasnikov_Gorbunov_2018}, where the authors proposed accelerated and non-accelerated zeroth-order method, which works in different proximal-setups. They obtained almost dimension-independent rate for the non-accelerated algorithm for the case of $\ell_1$-proximal-setup and sparse vector $x_0 - x_*$.

\subsection{Outline}

We organize this paper as follows. In Section~\ref{sec:RDS} we present our stochastic three points method and give some of its properties. In Section~\ref{sec:SD} we give the main assumptions on the random direction to ensure the convergence of our method. Then, in Section~\ref{sec:kl} we introduce the key lemma for the analysis of the complexity. Section~\ref{sec:ncc}
 gives the analysis for the worst case complexity for non-convex problems. While Section~\ref{sec:cc} deals with the complexity analysis for the convex problems, and Section~\ref{sec:scc} gives the analysis of the complexity for strongly convex problems. Section~\ref{sec:pa} proposes a parallel version of \texttt{STP} and gives the corresponding complexity analysis. 
 Numerical tests are illustrated and discussed in Section~\ref{sec:exper}. Conclusions and future improvements are discussed in
Section~\ref{sec:conc}.

\subsection{Notation}

Throughout this paper ${\cal D}$ will denote a probability distribution over $\R^n$. We use
 $\Exp \left[\cdot \right]$ to denote the expectation and $\ve{x}{y} = x^\top y$ corresponds to the inner product of $x$ and $y$. We denote also by 
 $\| \cdot \|_2$ the $\ell_2$-norm, and by $\| \cdot \|_{\cal D}$ a norm dependent on ${\cal D}$ which we introduce in Section~\ref{sec:SD}.


\section{Summary of contributions}

Here we highlight some of the contributions of this work.


{\bf A simple and flexible algorithm.} We study a novel variant of direct search based on random directions, which we call Stochastic Three Points (\texttt{STP}). 
It depends on at most  three  parameters: The starting point $x_0$ for the iterate, the probability distribution ${\cal D}$ on $\R^n$ to sample the directions, and in some cases an $\alpha_0$ 
to define the stepsize. 
 The probability distribution  ${\cal D}$ may be iteration dependent as far as it satisfies the required assumption (see Assumption~\ref{ass:P_sc}).  
In fact, Assumption~\ref{ass:P_sc} may be weakened by letting the probability distribution to depend on the iteration $k$  in the following way 
\begin{enumerate}
\item The quantity $\gamma_{{\cal D}_k} \eqdef \Exp_{s\sim {\cal D}_k} \|s\|_2^2$ is positive and uniformly bounded away from infinite.
 
\item There is a constant $\mu_{\cal D}>0$ and norm $\|\cdot\|_{\cal D}$ (independent from $k$) on $\R^n$ such that
\begin{equation}\label{eq:shs7hsk}\Exp_{s\sim {\cal D}_k}\; |\ve{g_k}{s}| \geq \mu_{\cal D} \|g_k\|_{\cal D},
\end{equation}
where $g_k = \nabla f(x_k)$.
\end{enumerate}
This assumption may be weakened even more by letting $\mu_{\cal D}$ and norm $\|\cdot\|_{\cal D}$ to dependent on $k$ and assuming i)
the uniform boundness of $\mu_{{\cal D}_k}$ away from zero, ii) and that $\|\cdot\|_{{\cal D}_k}$ is uniformly equivalent to a norm independent from $k$. To avoid unnecessary notations and for the sake of clarity and simplicity of the presentation, for the analysis we choose the probability distribution to be iteration independent  in this paper.
 
{\bf A general setting.}
Our approach covers some rather exotic first order methods:
\begin{itemize}
\item Normalized Gradient Descent (NGD) method: at iteration $k$, $s \sim {\cal D}$ means that 
$s = \frac{g_k}{\|g_k\|_2}$ with probability 1.

\item Signed Gradient Descent (SignGD) method: at iteration $k$, $s \sim {\cal D}$ means that $s = sign \left( g_k \right)$ with probability $1$, where the $sign$ operation is element wise sign.


\item  Normalized Randomized Coordinate Descent (NRCD) method (equivalently this method can be called also Randomized Signed Gradient Descent):  at iteration $k$, $s \sim {\cal D}$ means that $s = \frac{g_k^{i}}{ |g_k^{i} |}e_i$ if $g_k^{i} \neq 0$ and $s=0$ otherwise, with probability $\frac{1}{n}$, where $g_k^{i}$ is the $i-th$ component of $g_k$.

\item Normalized Stochastic Gradient Descent (NSGD) method: at iteration $k$, $s \sim {\cal D}$ means that 
$s = \hat{g}_k $ where $\hat{g}_k$ is  the  stochastic gradient satisfying $\Exp\left[ \hat{g}_k \right] = \frac{g_k}{\|g_k\|_2}$, and $\Exp\left[ \|\hat{g}_k\|_2^2  \right]  \le \sigma < \infty $.
\end{itemize}
The required assumption on ${\cal D}$ is satisfied in these cases (see Appendix~\ref{app:B}).

The probability distribution is also allowed to be either continuous or discrete,  so that we cover many known 
 strategies of choosing the directions in the DFO setting in the literature. For instance, if ${\cal D}$ is the uniform law on the unit sphere we recover the directions proposed in \cite{Karmanov_1974a, Karmanov_1974b, Polyak_1987, Nesterov_Spokoiny_2017}. If it is the discrete law on $\{ \pm e_i, i=1,\ldots,n\}$ we recover the directions proposed in \cite{Stich_2011}. 
 If it is the discrete law on $\{ \pm d_i, i=1,\ldots,n\}$ where $d_i,~i=1,\ldots,n$ form a basis of $\R^n$, \texttt{STP} can be seen as a random variant of the Simplified Direct Search (SDS) method studied in \cite{KR-DFO2014}.  

One of the main goals of flexibility in choosing the probability distribution ${\cal D}$ is the   efficiency for solving some  optimization problems which may have some specific properties like:
\begin{itemize}
\item The size of the problem to optimize is very large such that even the addition of two vectors may be unfeasible.
For instance if the dimension of the problem (i.e., the size of $x$) is larger than the available memory, then updating all the components of $x$ at each iteration is impossible. One is allowed to update only some components of $x$ at each iteration.

\item The objective function is not entirely defined at the beginning of the optimization process, like in the streaming optimization. In other words the data describing the objective function arrives in real time during the optimization process. At a given iteration (time)  
we can not evaluate the objective function in all points of $\R^n$. We can only 
evaluate the objective function in a set of directions (only some components of  $x$ can be updated).
\item Even if we have the entire objective function at the beginning of the optimization process, for some problems the  computation of the function value increases with the number of the perturbed variables. In other words, 
when perturbing all the components of $x$ the evaluation of $f$ takes a lot of time.
However by perturbing only one parameter (or a set of parameters) the objective is evaluated in reasonable time.
\item Some prior knowledge about Lipschitz constants in some directions is available. 
\end{itemize}
For these kind of situations the choices of ${\cal D}$ to be a continuous law is prohibited. 
However the discrete choices of ${\cal D}$ are the most convenient in these cases.

%

{\bf Practicality.} \texttt{STP} method is extremely simple to use in practice and its analysis is also simple compared to the state-of-the-art  direct search methods based on random directions/stepsizes. In fact, the most related work to \texttt{STP} is the  work in \cite{Nesterov_Spokoiny_2017}. In the latter work, the proposed stepsizes depend on the Lipschitz constant of the gradient of the objective function, which may not be known in practice. However, for \texttt{STP} we proposed several stepsize selection schemes. Some of them can be easily computed in practice.
 Moreover, our preliminary numerical experiments show that our approach is competitive in practice.

{\bf Better bounds.} We obtained compact worst case complexity bounds.  These bounds are similar to those obtained in \cite{Nesterov_Spokoiny_2017}. They depend 
linearly on the dimension of the considered problem, while this dependence is quadratic for deterministic direct search methods \cite{Vicente_2013, Vicente_2016, KR-DFO2014}.
In Table~\ref{tab:sumcompl} we summarize selected
complexity results (bounds on the number of function evaluations) obtained in this paper for \texttt{STP} method. 
In all cases we assume that $f$ is differentiable,
bounded below (by $f_*$), with $L$-Lipschitz gradient. The assumptions listed in the first column of the table are additional to this.
The quantity $R_0$ measures the size of a specific level set of $f$.  
The symbol $\propto $ means proportional. In fact, this symbol appears in the 
definition of the stepsizes, for instance $\alpha_k \propto \tfrac{1}{\sqrt{k+1}}$ means that $\alpha_k$ is equal to some constant $\alpha_0$ (independent from $k$) multiplied by $\tfrac{1}{\sqrt{k+1}}$. This constant $\alpha_0$ usually depends in the constants of the problem, like the Lipschitz constant and $x_0$.
More details about the 
definitions of all these quantities are given in the main text.

{\tiny{
\begin{table}[h]
\begin{center}
{\footnotesize
\begin{tabular}{|c|c|c|c|} 
  \hline
 \begin{tabular}{c} Assumptions on $f$\\(additional to \\
 $L$-smoothness)
 \end{tabular} & Stepsizes & Complexity  & Theorems \\
  \hline
    \hline
   & && \\
  none &
\begin{tabular}{c}
$\alpha_k \propto \tfrac{1}{\sqrt{k+1}}$ \\
 $\alpha_k \propto \e$
\end{tabular}
  &$O\left(\frac{n}{\e^2}\right)$& \ref{thm:nonconvex1},  \ref{thm:nonconvex2}  \\
   & && \\  
   \hline
   & && \\   
  \begin{tabular}{c}convex, \\ $R_0$ finite \end{tabular}&  
\begin{tabular}{c}
$\alpha_k \propto f(x_k) - f(x_*)$ \\
$\alpha_k \propto \frac{|f(x_{k}+ts_k)-f(x_k)|}{t}$
\end{tabular}    
    &$O\left(\frac{n}{\e}\right)$& \ref{thm:convex2}, \ref{thm:convex1}  \\
   & && \\    
  \hline
     & && \\
    \begin{tabular}{c}$\lambda$-strongly\\
     convex \end{tabular}&     
     \begin{tabular}{c}
 $\alpha_k \propto { (f(x_k) - f(x_*) )}^{\frac{1}{2}}$ \\
$\alpha_k \propto \frac{|f(x_k + ts_k) - f(x_k)|}{t}$
\end{tabular}  
     &$O\left(n\log\left(\frac{1}{\e}\right)\right)$ & \ref{thm:stronglyconvex2} , \ref{thm:stronglyconvex1}\\
     & && \\
  \hline
\end{tabular}
}
\end{center}
\caption{Summary of the complexity results obtained in this paper for \texttt{STP} method. Column ``Complexity'' defines the number of iterations needed to guarantee $\min_{k} \Exp \left[ \|\nabla f(x_k)\|_{\cal D} \right] \leq \e$ (second row) or $ \Exp \left[f(x_k) - f(x_*) \right]  \le \e$ (third and fourth rows).}
\label{tab:sumcompl}
\end{table}
}}

{\bf Parallel method.}
In Table~\ref{tab:sumcomplparallel} we summarize selected
complexity results (bounds on the number of function evaluations) obtained in this paper for the parallel version of the \texttt{STP} method.  
More details about the  definitions of all quantities appearing in the table are given in the main text. \texttt{PSTP} method gives the same rate as \texttt{STP} method with spherical setup but for wider range of distributions.

\begin{table}[h]
\begin{center}
{\footnotesize
\begin{tabular}{|c|c|c|c|} 
  \hline
 \begin{tabular}{c} Assumptions on $f$\\(additional to \\
 $L$-smoothness)
 \end{tabular} & Stepsizes & Complexity  & Theorems \\
  \hline
    \hline
   & && \\
  none &
\begin{tabular}{c}
$\alpha_k \propto \tfrac{1}{\sqrt{k+1}}$
\end{tabular}
  &$O\left(\frac{n}{\e^2}\right)$& \ref{thm:nonconvex_paral}  \\
   & && \\  
   \hline
   & && \\   
  \begin{tabular}{c}convex, \\ $R_0$ finite \end{tabular} 
    & 
\begin{tabular}{c}
$\alpha_k \propto f(x_k) - f(x_*)$ 
\end{tabular}    
    &$O\left(\frac{n}{\e}\right)$& \ref{thm:convex_paral}\\
   & && \\    
  \hline
     & && \\
    \begin{tabular}{c}$\lambda$-strongly\\
     convex \end{tabular} &     
     \begin{tabular}{c}
 $\alpha_k \propto { (f(x_k) - f(x_*) )}^{\frac{1}{2}}$ 
\end{tabular}  
     &$O\left(n\log\left(\frac{1}{\e}\right)\right)$ & \ref{thm:stronglyconvexparal} \\
     & &&\\
  \hline
\end{tabular}
}
\end{center}
\caption{Summary of the complexity results obtained in this paper for the parallel version of \texttt{STP} method. As before, column ``Complexity'' defines the number of iterations needed to guarantee $\min_{k} \Exp \left[ \|\nabla f(x_k)\|_{\cal D} \right] \leq \e$ (second row) or $ \Exp \left[f(x_k) - f(x_*) \right]  \le \e$ (third and fourth rows).}
\label{tab:sumcomplparallel}
\end{table}



{\bf Experiments.} 
We provide a number of experimental results, showing that our approach is a competitive algorithm in practice. In fact, we compared on a large set of problems our approach with the method (\ref{eq:Polyak-method}) as well as with the coordinate search method (the \texttt{DDS} method which uses the $2n$ coordinate directions). The experiments show that the use of the random directions 
 leads to a significant improvement in terms of the number of function evaluation. Indeed, our approach and method (\ref{eq:Polyak-method}) outperform the \texttt{DDS} method. Moreover, our approach 
 exhibits better performances than the other two methods. See Section~\ref{sec:exper} for a complete view on the experimental results.

\section{Stochastic Three Points method}
\label{sec:RDS}


Our {\em stochastic three points} (\texttt{STP}) algorithm is formalized below as Algorithm~\ref{alg:STP}.

\begin{algorithm}
\caption{\bf Stochastic Three Points (\texttt{STP})}
\vspace{-4ex}
\label{alg:STP}
\begin{rm}
\begin{description}
\item[]
\vspace{3ex}
\item[Initialization] \ \\
Choose $x_0\in \R^n$, stepsizes $\alpha_k>0$, probability distribution ${\cal D}$ on $\R^n$.
\vspace{1ex}
\item[For $k=0,1,2,\ldots$] \ \\
\vspace{-2ex}
\begin{enumerate}
\item Generate a random vector $s_k\sim {\cal D}$
\item Let $x_+ = x_k+\alpha_k s_k$ and $x_- = x_k - \alpha_k s_k$
\item  $x_{k+1} = \arg \min \{f(x_-), f(x_+),f(x_k)\}$
\end{enumerate}
\end{description}
\end{rm}
\end{algorithm}
%
%

Due to the randomness of the search directions $s_k$ for $k\ge 0$, the iterates 
 are also random vectors for all $k\ge 1$. The starting point $x_0$ is not random (the initial objective function value $f(x_0)$ is deterministic).
Note that \texttt{STP} never moves to a point with a larger objective value. This monotonicity property does not depend on ${\cal D}$ or the properties of $f$. Let us formulate this simple observation as a lemma.

\begin{lem}[Monotonicity]  \texttt{STP} produces a monotonic sequence of iterates, i.e., $f(x_{k+1})\leq f(x_k)$ for all $k\geq 0$. As a consequence, \begin{equation}\label{eq:monotonicity}\Exp[f(x_{k+1})\;|\; x_k] \leq f(x_k).\end{equation}
\end{lem}

Throughout the paper, we assume that $f$ is differentiable, bounded below and has $L$-Lipschitz gradient. 

\begin{ass}\label{ass:L-smooth} 
The objective function $f$ is $L$-smooth with $L>0$
 and bounded from below by $f_*\in \R$. 
  That is, $f$ has a Lipschitz continuous gradient with a Lipschitz constant $L$:
  $$   
\|\nabla f(x) - \nabla f(y) \|_2 \le L  \|x - y \|_2, \qquad \forall x,y \in \R^n 
  $$
  and $ f(x) \ge f_*$ for all $ x \in \R^n.$
\end{ass} 

\subsection{Random Search Directions}
\label{sec:SD}

Our analysis in the sequel of the paper will be based on the following key assumption.

\begin{ass}\label{ass:P} The probability distribution ${\cal D}$ on $\R^n$ has the following properties: 
\begin{enumerate}
\item The quantity $\gamma_{\cal D} \eqdef \Exp_{s\sim {\cal D}} \; \|s\|_2^2$ is positive and finite.
 
\item There is a constant $\mu_{\cal D}>0$ and norm $\|\cdot\|_{\cal D}$ on $\R^n$ such for all $g\in \R^n$,
\begin{equation}\label{eq:shs7hs}\Exp_{s\sim {\cal D}}\; |\ve{g}{s}| \geq \mu_{\cal D} \|g\|_{\cal D}.\end{equation}
\end{enumerate}
\end{ass}

Note that since all norms  in $\R^n$ are equivalent, the second part of the above assumption is satisfied if and only if
\[\inf_{\|g\|_2=1} \Exp_{s\sim {\cal D}}\; |\ve{g}{s}| >0.\]
However, as the next lemma illustrates, it will be convenient to work with norms that are allowed to depend on ${\cal D}$. We now give some examples of distributions for which the above assumption is satisfied.

\begin{lem} \label{lem1} Let $g\in \R^n$.  
\begin{enumerate}
\item If ${\cal D}$ is the uniform distribution on the unit sphere in $\R^n$, then
\begin{equation}
\gamma_{\cal D} = 1 \quad  \text{and} \quad \Exp_{s \sim {\cal D}}\;  | \<g, s> | \sim \frac{1}{\sqrt{2\pi n}} \|g\|_2.
\end{equation}
Hence, ${\cal D}$ satisfies Assumption~\ref{ass:P}  with $\gamma_{\cal D}=1$, $\|\cdot\|_{\cal D}=\|\cdot\|_2$ and $\mu_{\cal D} \sim \frac{1}{\sqrt{2\pi n}}$.

\item 

If ${\cal D}$ is the normal distribution with zero mean and identity over ${n}$ as covariance matrix. i.e., $s\sim N(0,\frac{I}{n})$,
then
\begin{equation}
\gamma_{\cal D} = 1 \quad \text{and} \quad \Exp_{s \sim {\cal D}}\;  | \<g, s> | = \frac{\sqrt{2}}{  \sqrt{n\pi}}\|g\|_2.
\end{equation}
Hence, ${\cal D}$ satisfies Assumption~\ref{ass:P}  with $\gamma_{\cal D}=1$, $\|\cdot\|_{\cal D}=\|\cdot\|_2$ and $\mu_{\cal D}= \frac{\sqrt{2}}{  \sqrt{n\pi}}$.

\item If ${\cal D}$ is the uniform distribution on  $\{e_1,\dots,e_n\}$, then
\begin{equation}
\gamma_{\cal D} = 1 \quad \text{and} \quad \Exp_{s\sim {\cal D}}\;  | \<g, s> |  = \frac{1}{n} \|g\|_1.
\end{equation}
Hence, ${\cal D}$ satisfies Assumption~\ref{ass:P}  with $\gamma_{\cal D}=1$, $\|\cdot\|_{\cal D}=\|\cdot\|_1$ and $\mu_{\cal D}=\tfrac{1}{n}$.
\item If ${\cal D}$ is an arbitrary distribution on  $\{e_1,\dots,e_n\}$ given by $P(s=e_i)=p_i  >0$, then
\begin{equation}
\gamma_{\cal D} = 1 \quad \text{and} \quad \Exp_{s\sim {\cal D}}\;  | \<g, s> |  = \|g\|_{\cal D} \eqdef \sum_{i=1}^n p_i |g_i|.
\end{equation}
Hence, ${\cal D}$ satisfies Assumption~\ref{ass:P}  with $\gamma_{\cal D}=1$  
 and $\mu_{\cal D}=1$.
\item 

If ${\cal D}$ is a distribution on $D= \{d_1,\ldots,d_n\}$ where $d_1,\ldots,d_n$ form an orthonormal basis of $\R^n$ and 
$P(s = d_i) = p_i$, then 
\begin{equation}
\gamma_{\cal D} = 1 \quad \text{and} \quad \Exp_{s\sim {\cal D}}\;  | \<g, s> |  = \|g\|_{\cal D} \eqdef \sum_{i=1}^n p_i |g_i|.
\end{equation}
Hence, ${\cal D}$ satisfies Assumption~\ref{ass:P}  with $\gamma_{\cal D}=1$  
 and $\mu_{\cal D}=1$.

\end{enumerate}

\end{lem}
\begin{proof}
See Appendix~\ref{app:A}.
\end{proof}


Without loss of generality, in the rest of this paper we assume that $\gamma_{\cal D}=1$. This can be achieved
 by considering distribution ${\cal D}'$  instead, where $s'\sim {\cal D}'$ is obtained by first sampling $s'$ from ${\cal D}$ and  then  either normalizing via i) $s = s'/ \|s'\|_2$, or ii) $s=s'/ \sqrt{ \Exp_{s'\sim {\cal D}} \|s'\|_2^2}$.
 
 \subsection{Key Lemma}
 \label{sec:kl}
Now, we establish the key result which will be used to prove the main properties of our Algorithm. Its similar result in the case of deterministic direct search (\texttt{DDS}) methods states that the gradient of the objective function for unsuccessful iterations is bounded by a constant multiplied by the stepsize. See for instance \cite[Lemma 10]{KR-DFO2014}.
\begin{lem}\label{lem:main} 

If Assumptions~\ref{ass:L-smooth} and~\ref{ass:P}  hold, then for all $k\geq 0$,
\begin{equation}\label{eq:s88ss}\Exp \left[ f(x_{k+1}) \;|\; x_k \right] \leq f(x_k)- \mu_{\cal D} \alpha_k \|\nabla f(x_k)\|_{\cal D} + \frac{L}{2}\alpha_k^2,
\end{equation}
and
\begin{equation}\label{eq:main-lemma}
\theta_{k+1} \leq \theta_k -\mu_{\cal D} \alpha_k g_k  + \frac{L}{2}\alpha_k^2,
\end{equation}
where  $\theta_k = \Exp [f(x_k)]$ and $g_k = \Exp [\|\nabla f(x_k)\|_{\cal D}]$.
\end{lem}
\begin{proof}
First we notice that from $L$-smoothness of $f$ we have
\begin{eqnarray*}
f(x_k+\alpha_k s_k) &\leq& f(x_k)+ \ve{\nabla f(x_k)}{\alpha_k s_k} + \tfrac{L}{2}\|\alpha_k s_k\|_2^2 \\
 &=& f(x_k) + \alpha_k \ve{\nabla f(x_k)}{s_k} + \tfrac{L}{2}\alpha_k^2\|s_k\|_2^2,
\end{eqnarray*}
and, similarly,
$f(x_k - \alpha_k s_k) \leq  f(x_k) - \alpha_k \ve{\nabla f(x_k)}{s_k} + \tfrac{L}{2}\alpha_k^2\|s_k\|_2^2$.
Hence,
\[f(x_{k+1}) \leq \min \{f(x_k+\alpha_k s_k) , f(x_k-\alpha_k s_k)\} \leq f(x_k) - \alpha_k |\ve{\nabla f(x_k)}{s_k}| + \tfrac{L}{2}\alpha_k^2\|s_k\|_2^2.\]
To conclude (\ref{eq:s88ss}), we only need to take expectation in the above inequality with respect to $s_k\sim {\cal D}$, conditional on $x_k$, and use inequality \eqref{eq:shs7hs}. By taking the expectation in (\ref{eq:s88ss}) we get (\ref{eq:main-lemma}).
\end{proof}

Note that \eqref{eq:s88ss} can equivalently be written in the following form:

\[\|\nabla f(x_k)\|_{\cal D} \leq \frac{1}{\mu_{\cal D} } \left( \frac{f(x_k) - \Exp \left[ f(x_{k+1}) \;|\; x_k \right]}{\alpha_k} +  \frac{L}{2}\alpha_k\right).\]

This form makes it possible to compare this result with a key result used in the analysis of \texttt{DDS}. Indeed, if we assume that the opposite of the following sufficient {\em expected} decrease condition holds
\begin{equation}
\label{eq:nsd}
f(x_k) - \Exp \left[ f(x_{k+1}) \;|\; x_k \right] \geq c\alpha_k^2,
\end{equation}
for some  $c>0$, then we obtain
\begin{equation}
\label{eq:ugrad}
\|\nabla f(x_k)\|_{\cal D} \leq \frac{1}{\mu_{\cal D} } \left( c +  \frac{L}{2}\right)\alpha_k.
\end{equation}
In \texttt{DDS}, condition (\ref{eq:nsd}) is equivalent to the sufficient decrease condition $f(x_k) -  f(x_{k+1}) \geq c\alpha_k^2$. If such condition does not hold than the step is declared unsuccessful.
The inequality in (\ref{eq:ugrad}) is similar with the result in \cite[Lemma 10]{KR-DFO2014}. In \texttt{DDS} methods, one can check the sufficient decrease condition, so this drives the analysis and allows for simple stepsize update rules to be implemented. In \texttt{STP}, we typically cannot evaluate $\Exp[f(x_{k+1})\;|\; x_k]$ (we can if ${\cal D}$ has all its mass on a discrete set -- but in that case we would need to do more work per iteration).

\section{Non-convex Problems}
\label{sec:ncc}

In this section, we state our most general complexity result where we do not make  any additional assumptions on $f$, besides smoothness and boundedness (see Assumption~\ref{ass:L-smooth}). 

\begin{thm}[Decreasing stepsize]\label{thm:nonconvex1} Let Assumptions~\ref{ass:L-smooth} and~\ref{ass:P} hold. Choose $\alpha_k=\tfrac{\alpha_0}{\sqrt{k+1}}$, where $\alpha_0>0$. If \begin{equation}\label{eq:isjss8sus} K\geq \frac{2\left(\frac{\sqrt{2}(f(x_0)-f_*)}{\alpha_0} + \frac{L\alpha_0}{2}\right)^2}{\mu_{\cal D}^2 \e^2},\end{equation}
then  $\min_{k=0,1,\dots,K} \Exp \left[ \|\nabla f(x_k)\|_{\cal D} \right] \leq \e.$
 
\end{thm}

\begin{proof}
We base the proof on the analysis of the recursion \eqref{eq:main-lemma}. In particular, it is useful to write it in the following form:
\begin{equation}\label{eq:s9jd7d76d} g_k \leq \tfrac{1}{\mu_{\cal D}} \left( \tfrac{\theta_k-\theta_{k+1}}{\alpha_k} + \tfrac{L}{2} \alpha_k \right) = \tfrac{1}{\mu_{\cal D}} \left( \tfrac{(\theta_k-\theta_{k+1})\sqrt{k+1}}{\alpha_0} + \tfrac{L\alpha_0}{2 \sqrt{k+1}}\right).\end{equation}
We know from \eqref{eq:monotonicity} and the assumption that $f$ is bounded below that $f_*\leq \theta_{k+1}\leq \theta_k\leq f(x_0)$ for all $k$. Letting $l=\lfloor K/2 \rfloor$, this implies that \[\sum_{j=l}^{2l} (\theta_{j}-\theta_{j+1})  = \theta_l - \theta_{2l+1} \leq  f(x_0)-f_*\eqdef C,\]
from which we conclude that there must exist $j\in \{l,\dots,2l\}$ such that $\theta_j-\theta_{j+1}\leq C/(l+1)$. This implies that
\begin{eqnarray*}
g_j &\overset{\eqref{eq:s9jd7d76d}}{\leq} & \tfrac{1}{\mu_{\cal D}} \left( \tfrac{(\theta_j-\theta_{j+1})\sqrt{j+1}}{\alpha_0} + \tfrac{L\alpha_0}{2 \sqrt{j+1}}\right) \leq  \tfrac{1}{\mu_{\cal D}} \left( \tfrac{C\sqrt{j+1}}{\alpha_0 (l+1)} + \tfrac{L\alpha_0}{2 \sqrt{j+1}}\right)\\
&\leq & \tfrac{1}{\mu_{\cal D}} \left( \tfrac{C\sqrt{2l+1}}{\alpha_0 (l+1)} + \tfrac{L\alpha_0}{2 \sqrt{l+1}}\right) \leq  \tfrac{1}{\mu_{\cal D} \sqrt{l+1}} \left( \tfrac{\sqrt{2}C}{\alpha_0 } + \tfrac{L\alpha_0}{2}\right)\\
&\leq & \tfrac{1}{\mu_{\cal D} \sqrt{K/2}} \left( \tfrac{\sqrt{2}C}{\alpha_0 } + \tfrac{L\alpha_0}{2}\right) \overset{\eqref{eq:isjss8sus}}{\leq} \e.
\end{eqnarray*}
\end{proof}

Let us now give some insights into the above theorem.
\begin{itemize}
\item \textbf{Sphere setup.} If ${\cal D}$ is the uniform distribution on the Euclidean sphere, then  $\mu_{\cal D} \sim \frac{1}{\sqrt{2\pi n}}$, and hence the above theorem gives a complexity guarantee of the form
\[O\left(\frac{n}{\e^2}\right).\]
This is an improvement on \texttt{DDS} where the best known complexity bound is $O(n^2/\e^2)$ \cite{Vicente_2013, KR-DFO2014}. The same conclusion holds for the normal distribution setup.

\item \textbf{Coordinate setup.} If ${\cal D}$ is the uniform distribution on $\{e_1,\dots,e_n\}$, then $\mu_{\cal D}=1/n$ and hence the bound is of the form
\[O\left(\frac{n^2}{\e^2}\right).\]
However, this is for the $\ell_1$ norm of the gradient of $f$, which is {\em larger} than the $\ell_2$ norm. Indeed, for all $x$ we have $\sqrt{n}\|\nabla f(x)\|_2 \geq \|\nabla f(x)\|_1 \geq \|\nabla f(x)\|_2$, and the first inequality can be tight (for the vector of all ones, for instance). Hence, if we are interested to achieve $\|\nabla f(x)\|_2\leq \e'$, in certain situations it may be sufficient to push the $\ell_1$ norm of the gradient below $\e=\sqrt{n}\e'$ instead. So, the iteration bound can be as good as
\[O\left(\frac{n^2}{(\sqrt{n}\e')^2}\right) = O\left(\frac{n}{(\e')^2}\right).\]

\item \textbf{Quality of the final iterate.} Theorem~\ref{thm:nonconvex1} does not guarantee the gradient of $f$ at  the {\em final} point $x_K$ to be small (in expectation). Instead, it guarantees that the gradient of $f$ at {\em some} point produced by the method will be small. Notice however, that the method is monotonic. Hence, all subsequent points produced by the method will have better functions values  than the one which has gradient of minimum norm (in expectation). So, we can say that 
$f(x_K) \leq f(x_j)$ where $\Exp \left [\|\nabla f(x_j)\|_{\cal D} \right]\leq \e$.
\item \textbf{Optimal stepsize.} Note that the complexity depends on $\alpha_0$. The optimal choice (minimizing the complexity bound) is 
\[\alpha^* = 8^{1/4}\sqrt{\frac{f(x_0)-f_*}{L}},\]
in which case the complexity bound \eqref{eq:isjss8sus} takes the form
\begin{equation}\label{eq:isjss8sus-optimal} \frac{4\sqrt{2}(f(x_0)-f_*)L}{\mu_{\cal D}^2 \e^2}.\end{equation}
Assume that the lower bound  $f_*$ is achieved by some point $x_*\in \R^n$. Necessarily, $\nabla f(x_*)=0$. Moreover, since $f$ is $L$-smooth, we can write
\[f(x_0)\leq f(x_*) + \ve{\nabla f(x_*)}{x_0-x_*} + \frac{L}{2}\|x_0-x_*\|_2^2.\]
Hence, the optimal stepsize is no larger than
\[\alpha^* \leq 2^{1/4} \|x_0-x_*\|_2.\]

 Of course, we cannot use this optimal stepsize as we usually do not know $L$ and/or $f_*$. So, we are paying for the lack of knowledge by an increased complexity bound. This makes intuitive sense: the stepsize should not be much larger than the distance of the initial point to an optimal point.

On the other hand, there are examples of non-convex functions for which  the ratio $(f(x_0)-f_*)/L$ is arbitrarily small, and the distance between $x_0$ and $x_*$ arbitrarily high. This cannot happen for convex functions with bounded level sets or for strongly convex functions, as then $f(x) -  f(x_*)$ can be lower bounded by quantity proportional to $\|x - x_*\|_2$ with some positive power.

 \end{itemize}

We now state a complexity theorem for \texttt{STP} used with a fixed stepsize.

\begin{thm}[Fixed stepsize]\label{thm:nonconvex2} Let $f$ satisfy Assumption~\ref{ass:L-smooth} and also assume that $f$ is bounded below by $f_*\in \R$.  Choose a fixed stepsize $\alpha_k=\alpha$ with $0<\alpha< 2 \mu_{\cal D} \e / L$. If \begin{equation}\label{eq:isjsus}K\geq k(\e)\eqdef  \left \lceil \frac{f(x_0)-f_*}{(\mu_{\cal D}\e - \tfrac{L}{2}\alpha)\alpha} \right \rceil -1,\end{equation}
then  $\min_{k=0,1,\dots,K} \Exp \left[ \|\nabla f(x_k)\|_{\cal D} \right] \leq \e.$ In particular, if $\alpha=\mu_{\cal D}\e/L $, then
\[k(\e) = \left\lceil\frac{2L(f(x_0)-f_*)}{\mu_{\cal D}^2 \e^2}\right\rceil -1.\]
\end{thm}
\begin{proof}  If  $g_k\le \e $ for some $k\leq k(\e)$, then we are done. Assume hence by contradiction that $g_k>\e$ for all $k\leq k(\e)$.  By taking expectation in  Lemma~\ref{lem:main}, we get
\[\theta_{k+1} \leq \theta_k -\mu_{\cal D} \alpha g_k +\tfrac{L}{2}\alpha^2,\]
where $\theta_k = \Exp [f(x_k)]$ and $g_k = \Exp [\|\nabla f(x_k)\|_{\cal D}]$. Hence, 
\[f_* \leq \theta_{K+1} < \theta_0 - (K+1) \left(\mu_{\cal D} \alpha \e - \tfrac{L}{2}\alpha^2\right) \overset{\eqref{eq:isjsus}}{\leq} \theta_0 - (f(x_0)-f_*) = f_*,\]
which is a contradiction.

\end{proof}

Here we give some comments about the \texttt{STP} for non-convex functions.

\begin{itemize}
\item  In some situations, when $L$ is not available, it is impossible to compute optimal $\alpha = \frac{\mu_{\cD}\e}{L}$.
\item  If we can guess $\alpha$ is close to the optimal, then the method depends linearly on $n$ if $1/\mu_{\cal D}^2 = O(n)$. 
\item Also, if we guess $\alpha$ right, we get complexity that depends on $L(f(x_0)-f_*)$, which is similar to the setup with variable stepsizes and optimal $\alpha_0$.
\item As before, we only get guarantee on the best of the points in term of the gradient norm, not on the final point.
\end{itemize}

\section{Convex Problems}
\label{sec:cc}

In this section we estimate the complexity of the \texttt{STP}  in the case of convex $f$. In this case we need an additional technical assumption.

\begin{ass} \label{ass:levelset} We assume that $f$ is convex, has a minimizer $x_*$, and has bounded level set at $x_0$: 
\[R_0 \eqdef \max \{ \|x-x_*\|_{\cal D}^*\;:\; f(x)\leq f(x_0)\}< +\infty,\]
where $\|\xi\|_{\cal D}^* \eqdef \max\{\langle\xi, x\rangle\mid \|x\|_{\cal D} \leq 1 \}$ defines the dual norm to $\|\cdot\|_{\cal D}$.
\end{ass}

Note that if the above assumption holds, then whenever $f(x)\leq f(x_0)$, we get $f(x)-f(x_*)\leq \ve{\nabla f(x)}{x-x_*}\leq \|\nabla f(x)\|_{\cal D}\|x-x_*\|_{\cal D}^* \leq R_0 \|\nabla f(x)\|_{\cal D}$. That is, 
\begin{equation}\label{eq:s8jd7dh}\|\nabla f(x)\|_{\cal D} \geq \frac{f(x)-f(x_*)}{R_0}.\end{equation}

Now, we state our main complexity result of this section. We start with the analysis of \texttt{STP} with constant stepsizes.
\begin{thm}[Constant stepsize]   \label{thm:convex3}
\label{thm:aasdfsafa}
Let Assumptions~\ref{ass:L-smooth},~\ref{ass:P}~and~\ref{ass:levelset} be satisfied. Let $0< \e< \frac{L R_0^2}{\mu_{\cal D}^2 }$ and choose constant stepsize $\alpha_k =\alpha= \tfrac{\e \mu_{\cal D} }{LR_0}$.  If 
\begin{equation}\label{eq:isjsssus}K\geq \frac{LR_0^2}{\mu_{\cal D}^2  \e} \log \left(\frac{2(f(x_0)-f(x_*))}{\e}\right),
\end{equation}
then  $ \Exp \left[f(x_K) - f(x_*) \right] \leq \e.$
\end{thm}
\begin{proof}  Let us substitute \eqref{eq:s8jd7dh} into  Lemma~\ref{lem:main} and take expectations. We get
\begin{equation}\label{eq:ss8sjs8}\theta_{k+1}  \leq \theta_k  - \tfrac{\mu_{\cal D} \alpha}{R_0} (\theta_k - f(x_*)) +\tfrac{L}{2}\alpha^2.\end{equation}
 Let $r_k = \theta_k-f(x_*)$ and $c = 1-\tfrac{\mu_{\cal D} \alpha}{R_0}\in (0,1)$. Subtracting $f(x_*)$ from both sides of \eqref{eq:ss8sjs8}, we obtain 
\begin{eqnarray*}r_K &\leq& c r_{K-1}  +\tfrac{L}{2}\alpha^2 \quad 
\leq \quad c^K r_0+ \tfrac{L}{2}\alpha^2\sum_{i=0}^{K-1} c^i \\
&\leq & \exp(-\mu_{\cal D} \alpha K/ R_0) r_0 + \tfrac{L\alpha^2}{2(1-c)} 
\quad = \quad \exp(-\mu_{\cal D} \alpha K/ R_0) r_0 + \tfrac{\e}{2} \;\; \overset{\eqref{eq:isjsssus}}{\leq} \;\; \e.
\end{eqnarray*}
\end{proof}

If   $\mu_{\cal D} \sim \frac{1}{\sqrt{ n}}$, then the above theorem gives a complexity guarantee of the form
\[O\left(\frac{n}{\e}\log\left(\frac{1}{\e}\right)\right).\]
Comparing this to the best known complexity bound for \texttt{DDS} which is $O(\frac{n^2}{\e})$ \cite{Vicente_2016, KR-DFO2014}, 
we improve the dependence on $n$ but we deteriorate the dependence on $\e$ because of the presence of the term $\log\left(\frac{1}{\e}\right)$.
 In the next theorem we show how we get rid of the $\log\frac1\e$ term using variable stepsize.

\begin{thm}[Variable stepsize] \label{thm:convex2} Let Assumptions~\ref{ass:L-smooth},~\ref{ass:P}~and~\ref{ass:levelset} be satisfied. 
Let $\alpha_k = \alpha_0 \left(f(x_k) - f(x_*)\right)$,  
 where $0 < \alpha_0 < \frac{2\mu_{\cal D}}{R_0 L}$.  Define $a  = \frac{\mu_{\cal D}\alpha_0}{R_0} - \frac{L \alpha_0^2}{2}>0.$
  If $k \ge k(\e) \eqdef \frac{1}{a}\left(\frac{1}{\e} - \frac{1}{r_0}\right),$ then $ \Exp \left[f(x_k) - f(x_*) \right]  \le \e$. 
 \end{thm}
 \begin{proof}
 Let us substitute \eqref{eq:s8jd7dh} into equation (\ref{eq:s88ss}) of  Lemma~\ref{lem:main}, and then substrate $f(x_*) $ from both sides  we get
\begin{equation*}
\Exp \left[ f(x_{k+1}) \;|\; x_k \right] - f(x_*) \leq f(x_k) -f(x_*)  - \mu_{\cal D} \alpha_k \tfrac{f(x_k)-f(x_*)}{R_0} + \tfrac{L}{2}\alpha_k^2.
\end{equation*} 
Let $r_k = \Exp \left[f(x_k)\right] - f(x_*) $. By using our choice of $\alpha_k$ in the previous equation and then taking the expectation we get
$
r_{k+1} \le r_k -\left( \tfrac{\mu_{\cal D}\alpha_0}{R_0} - \tfrac{L \alpha_0^2}{2}\right) r_k^2 = r_k - a r_k^2.
$ Therefore,
\begin{align*}
\tfrac1{r_{k+1}} 
-\tfrac1{r_k}
 =
 \tfrac{r_k - r_{k+1}}{r_k r_{k+1} }
 \geq 
 \tfrac{r_k - r_{k+1}}{r_k^2 } 
\geq a.
\end{align*} 
From this
we have
$
\tfrac1{r_k} \geq  \tfrac1{r_0} + ka
$
and hence
$
r_k \le \tfrac1{\tfrac1{r_0} + ka}.
$
It remains to notice that for $k \ge \tfrac{1}{a}\left(\tfrac{1}{\e} - \tfrac{1}{r_0}\right)$ we have 
$r_k  \le  \tfrac1{\tfrac1{r_0} + ka} \le  \e$.
 \end{proof}
 If $\alpha_0  = \frac{\mu_{\cal D}}{R_0L}$, then $a$ is maximal as a function of $\alpha_0$, for which we get the optimal bound  $$k(\e)=  \frac{2R_0^2L}{\mu_{\cal D}^2}\left(\frac{1}{\e} - \frac{1}{r_0}\right).$$
If   $\mu_{\cal D} \sim \frac{1}{\sqrt{ n}}$, then the above theorem gives a complexity guarantee of the form $O\left(\frac{n}{\e}\right).$

 The stepsizes in the previous theorem depend on $f(x_*)$. Of course, in practice we cannot always use these stepsizes as we usually do not know $f(x_*)$. Next theorem gives a more practical stepsizes for which we get the same complexity as in the previous theorem. 
We start by stating an extra assumption on the probability distribution ${\cal D}$ and show that this assumption is satisfied for all the probability distributions  given in  Lemma~\ref{lem1}.

\begin{ass}\label{ass:P_sc} The probability distribution ${\cal D}$ on $\R^n$ is such that for all $s\sim {\cal D}$ are of unit Euclidean norm $(\|s\|_2=1)$ with probability 1.
\end{ass}
Let $C_{\cal D}$ be the positive constant such that for all $x\in\R^n$ the following inequality holds: $\|x\|_{2} \leq C_{\cal D}\|x\|_{\cD}$. Such constant exists due to the equivalence of the norms in $\R^n$. 

\begin{thm}[Solution-free stepsize]\label{thm:convex1} Let Assumptions~\ref{ass:L-smooth},~\ref{ass:P},~\ref{ass:levelset}~and~\ref{ass:P_sc} be satisfied. Let $\alpha_k = \frac{|f(x_{k}+ts_k)-f(x_k)|}{Lt}$, where \[0< t \le \frac{\sqrt{2}\mu_{\cD} \Exp\left[ f(x_{K-1}) - f_* \right]}{L R_0}.\] Define $a  = \frac{\mu_{\cal D}^2}{4LR_0^2}.$ If $K \ge k(\e) \eqdef \frac{1}{a}\left(\frac{1}{\e} - \frac{1}{r_0}\right),$ then $ \Exp \left[f(x_K) - f(x_*) \right]  \le \e$. 
\end{thm}
\begin{proof}
	 From Lemma~\ref{lem:main} we have
	\begin{equation}\label{eq:lsmth_func_eval_conv}
	f(x_{k+1}) \le f(x_k) - \alpha_k|\langle \nabla f(x_k), s_k \rangle| + \tfrac{L\alpha_k^2}{2}.
	\end{equation}
	We know that $\alpha_k^{\text{opt}} = \tfrac{|\langle \nabla f(x_k), s_k \rangle|}{L}$ minimizes the right hand side of \eqref{eq:lsmth_func_eval_conv}. But it depends on $\nabla f(x_k)$ which we can not compute \textit{exactly}, because we have zeroth-order oracle. Actually, we do not need to know the whole gradient, it is enough to know directional derivative of $f$, which we can approximate by finite difference of function values of $f$. It is the main idea behind our choice of $\alpha_k^{\text{opt}} = \tfrac{|f(x_k+ts_k) - f(x_k)|}{Lt},$ which does not depends any more on $f(x_*)$ and can be easily  computed in practice. We can rewrite $\alpha_k = \tfrac{|f(x_k + ts_k) - f(x_k)|}{L t} = \tfrac{|\langle\nabla f(x_k),s_k\rangle|}{L} + \tfrac{|f(x_k + ts_k) - f(x_k)|}{L t} - \tfrac{|\langle\nabla f(x_k),s_k\rangle|}{L} \eqdef \alpha_k^{\text{opt}} + \delta_k$. Therefore, we have
	\begin{eqnarray*}
		f(x_{k+1}) &\le& f(x_k) - \tfrac{|\langle\nabla f(x_k),s_k\rangle|^2}{L} - \delta_k|\langle\nabla f(x_k),x_k\rangle| + \tfrac{|\langle\nabla f(x_k),x_k\rangle|^2}{2L}\\
		&&\quad + \delta_k|\langle\nabla f(x_k),x_k\rangle| + \tfrac{L}{2}(\delta_k)^2\\
		&=& f(x_k) - \tfrac{|\langle\nabla f(x_k),s_k\rangle|^2}{2L} + \tfrac{L}{2}(\delta_k)^2
	\end{eqnarray*}
	Next we estimate $|\delta_k|$ using $L$-smoothness of $f$:
	\begin{eqnarray*}
		|\delta_k| &=& \tfrac{1}{L t}\left||f(x_k + ts_k) - f(x_k)| - |\langle\nabla f(x_k),ts_k\rangle|\right|\\
		&\le& \tfrac{1}{L t}\left|f(x_k + ts_k) - f(x_k) - \langle\nabla f(x_k),ts_k\rangle\right|
		\le \tfrac{1}{L t} \cdot\tfrac{L}{2}\|ts_k\|_2^2 = \tfrac{t}{2}.
	\end{eqnarray*}		
	From this we obtain
	\begin{equation}
	\label{eq:convexos}
	\begin{array}{cl}
	f(x_{k+1}) &\leq f(x_k) - \tfrac{\left|\langle\nabla f(x_k), s_k \rangle\right|^2}{2L} + \tfrac{Lt^2}{8}.
	\end{array}
	\end{equation}
	Taking mathematical expectation w.r.t. all randomness from the previous inequality we get
	\begin{equation}\label{eq:conv_main_estim}
	\begin{array}{cl}
	\underbrace{\Exp[f(x_{k+1})] - f_*}_{r_{k+1}} & \overset{\circledOne}{\leq} \underbrace{\Exp[f(x_k)] - f_*}_{r_k} - \tfrac{\mu_{\cal D}^2}{2L}\Exp[\|\nabla f(x_k)\|_{\cal D}^2] + \tfrac{Lt^2}{8}\\
	&\overset{\circledTwo}{\leq} r_k - \tfrac{\mu_{\cal D}^2}{2LR_0^2}r_k^2 + \tfrac{Lt^2}{8},
	\end{array}
	\end{equation}
	where $\circledOne$ is due to tower property of mathematical expectation and \eqref{eq:shs7hs}:
	\begin{eqnarray*}
		\Exp[|\langle \nabla f(x_k), s_k \rangle|^2] &=& \Exp\left[\Exp[|\langle \nabla f(x_k), s_k \rangle|^2 \mid x_k]\right] \geq \Exp\left[\left(\Exp[|\langle \nabla f(x_k), s_k \rangle| \mid x_k]\right)^2\right]\\
		&\overset{\eqref{eq:shs7hs}}{\geq}& \mu_{\cal D}^2\Exp[\|\nabla f(x_k)\|_{\cal D}^2];	
	\end{eqnarray*}
	$\circledTwo$ follows from Assumption~\ref{ass:levelset}: $\Exp[\|\nabla f(x_{k})\|_{\cal D}^2] \geq \tfrac{\Exp\left[\left(f(x_k)-f_*\right)^2\right]}{R_0^2} \geq \tfrac{\left(\Exp\left[f(x_k) - f_*\right]\right)^2}{R_0^2} = \tfrac{r_k^2}{R_0^2}$. From this and monotonicity of $\{f(x_k)\}_{k\ge 0}$ we have
	\begin{equation}\label{eq:conv_pre_final}
		\tfrac{1}{r_{k+1}} - \tfrac{1}{r_k} \geq \tfrac{r_{k+1}-r_k}{r_k r_{k+1}} \geq \tfrac{\tfrac{\mu_{\cD}^2}{2LR_0^2}r_k^2 - \tfrac{Lt^2}{8}}{r_k^2} \geq \tfrac{\mu_{\cD}^2}{2LR_0^2} - \tfrac{L}{8}\left(\tfrac{t}{r_k}\right)^2.
	\end{equation}
	If $k \le K-1$ and $0< t \leq \tfrac{\sqrt{2}\mu_{\cD}r_{K-1}}{LR_0},$ then we can write
	\begin{eqnarray*}
		\tfrac{1}{r_{k+1}} - \tfrac{1}{r_k} \ge \tfrac{\mu_{\cD}^2}{4LR_0^2} = a,
	\end{eqnarray*}
	since $r_k \le r_{K-1}$.
	Finally, we have
$
	\tfrac1{r_k} \geq  \tfrac1{r_0} + ka
$
	and hence
	$
	r_k \le \tfrac1{\tfrac1{r_0} + ka}.
	$ for all $k \le K$.
	Thus, if $K \ge \tfrac{1}{a}\left(\tfrac{1}{\e} - \tfrac{1}{r_0}\right),$ then 
	$r_K  \le  \tfrac1{\tfrac1{r_0} + Ka} \le  \e$.
\end{proof}

Actually, requirement $t \le \frac{\sqrt{2}\mu_{\cD} \Exp\left[ f(x_{K-1}) - f_* \right]}{L R_0}$ could be replaced by $t \le \frac{\sqrt{2}\mu_{\cD} \e}{L R_0}$ if we additionally require that for all $k \le K$ we have $r_k \ge \e$.

\section{Strongly Convex Problems}
\label{sec:scc}

In this section we derive the complexity of the \texttt{STP} method in the case of strongly convex $f$.  

\begin{ass}\label{ass:strongconvex} $f$ is $\lambda$-strongly convex with respect to the norm $\|\cdot\|_{\cal D}$. 
\end{ass}

In this section, we denote by  $x_*$ the unique  minimizer of $f$.

\begin{thm}  \label{thm:stronglyconvex2}
Let Assumptions~\ref{ass:L-smooth},~\ref{ass:P}~and~\ref{ass:strongconvex} be satisfied. Let stepsize $\alpha_k = \frac{\theta_k \mu_{\cal D}}{L}\sqrt{2\lambda (f(x_k) - f(x_*) )}$, 
 for some $\theta_k \in (0,2)$ such that $\theta\eqdef \inf_k 2\theta_k-\theta_k^2>0$.
 If    \begin{equation}\label{eq:isjsssusc}K\geq \frac{L}{\lambda \mu_{\cal D}^2  \theta} \log \left(\frac{f(x_0)-f(x_*)}{\e}\right),\end{equation}
  then $\Exp \left[f(x_K) - f(x_*) \right]  \le \e$.
\end{thm}
\begin{proof}

By injecting $\alpha_k$ into equation (\ref{eq:s88ss}) of  Lemma~\ref{lem:main}, and then substrate $f(x_*) $ from both sides  we get
\begin{eqnarray*}
\Exp[f(x_{k+1}) \;|\; x_k] -f(x_*)  &\le& f(x_k) - f(x_*)  -\tfrac{\mu_{\cal D}^2\theta_k \sqrt{2\lambda( f(x_k) - f(x_*) )}\|\nabla f(x_k) \|_{\cal D}}{L}\\
&& \quad + \tfrac{\mu_{\cal D}^2\theta_k^2 \lambda (f(x_k) - f(x_*) )}{L}.
\end{eqnarray*}
From the strong convexity property of $f$ we have $\|\nabla f(x_k) \|_{\cal D}^2 \ge 2\lambda (f(x_k) - f(x_*) ),$ therefore
\begin{eqnarray*}  
\Exp[f(x_{k+1}) | x_k ] -f(x_*)  &\le& f(x_k) - f(x_*)  -\tfrac{2\mu_{\cal D}^2\theta_k \lambda( f(x_k) - f(x_*) )}{L} + \tfrac{\mu_{\cal D}^2\theta_k^2 \lambda (f(x_k) - f(x_*) )}{L}\\
&\le & f(x_k) - f(x_*)  -\tfrac{\mu_{\cal D}^2 \lambda( f(x_k) - f(x_*) )}{L} (2\theta_k - \theta_k^2)\\
&\le & f(x_k) - f(x_*)  -\tfrac{\mu_{\cal D}^2\theta \lambda( f(x_k) - f(x_*) )}{L}, 
\end{eqnarray*}
where we used the definition of $ \theta$. Let $r_k = \Exp \left[f(x_k)\right] - f(x_*) $. By taking the expectation of the last inequality we get
$
r_{k+1}   \le \left(1 - \tfrac{\mu_{\cal D}^2\theta\lambda}{L}\right) r_k,
$
and therefore 
\begin{eqnarray*}
r_{k}   \le \left(1 - \tfrac{\mu_{\cal D}^2\theta\lambda}{L}\right)^k r_0. 
\end{eqnarray*}
Hence if $K$ satisfies (\ref{eq:isjsssusc}), we get $r_K \le \e$.
\end{proof}

From this theorem we conclude that if there exist $0<\theta_1 \le \theta_2<2$ such that
 $$ \frac{\theta_1 \mu_{\cal D}}{L}\sqrt{2\lambda (f(x_k) - f(x_*) )} \le \alpha_k \le \frac{\theta_2 \mu_{\cal D}}{L}\sqrt{2\lambda (f(x_k) - f(x_*) )},$$
 then the sequence $(r_k)_k$ converges linearly to zero.
 
 The stepsizes from the previous theorem depend on $f(x_*)$. In practice,  we cannot always use these stepsizes as we usually do not know $f(x_*)$. Next theorem gives the similar result for \texttt{STP} with stepsizes independent from $f(x_*)$ under additional assumptions that for all $s\sim\cD$ we have $\|s\|_2 = 1$ with probability $1$.

\begin{thm}   \label{thm:stronglyconvex1}
	Let Assumptions~\ref{ass:L-smooth},~\ref{ass:P},~\ref{ass:P_sc}~and~\ref{ass:strongconvex} be satisfied.
	Let $\alpha_k = \frac{|f(x_k + ts_k) - f(x_k)|}{Lt}$, 
	for
$
		1 < t \le \frac{2\mu_{\cD}\sqrt{\lambda\e}}{L}.
$
	If    \begin{equation}\label{eq:num_it}K\geq \frac{L}{\lambda \mu_{\cal D}^2} \log \left(\frac{2(f(x_0)-f(x_*))}{\e}\right),\end{equation}
	then $\Exp \left[f(x_K) - f(x_*) \right]  \le \e$.
\end{thm}
\begin{proof}
From (\ref{eq:convexos}) we have
$
	f(x_{k+1}) 
	\leq f(x_k) - \tfrac{\left|\langle\nabla f(x_k), s_k \rangle\right|^2}{2L} + \tfrac{Lt^2}{8}.
$
	Taking mathematical expectation w.r.t. all randomness from the previous inequality we get
	\begin{equation}\label{eq:sc_main_estim}
	\begin{array}{cl}
	\underbrace{\Exp[f(x_{k+1})] - f_*}_{r_{k+1}} & \overset{\circledOne}{\leq} \underbrace{\Exp[f(x_k)] - f_*}_{r_k} - \tfrac{\mu_{\cal D}^2}{2L}\Exp[\|\nabla f(x_k)\|_{\cal D}^2] + \tfrac{Lt^2}{8}\\
	&\overset{\circledTwo}{\leq} \left(1-\tfrac{\mu_{\cal D}^2\lambda}{L}\right)r_k + \tfrac{Lt^2}{8},
	\end{array}
	\end{equation}
	where $\circledOne$ is due to tower property of mathematical expectation and \eqref{eq:shs7hs}:
	\begin{eqnarray*}
		\Exp[|\langle \nabla f(x_k), s_k \rangle|^2] &=& \Exp\left[\Exp[|\langle \nabla f(x_k), s_k \rangle|^2 \mid x_k]\right] \geq \Exp\left[\left(\Exp[|\langle \nabla f(x_k), s_k \rangle| \mid x_k]\right)^2\right]\\
		&\overset{\eqref{eq:shs7hs}}{\geq}& \mu_{\cal D}^2\Exp[\|\nabla f(x_k)\|_{\cal D}^2];	
	\end{eqnarray*}
	$\circledTwo$ follows from $\lambda$-strong convexity of $f$: $\|\nabla f(x_k)\|_{\cD}^2 \ge 2\lambda\left(f(x_k)-f_*\right)$. From \eqref{eq:sc_main_estim} we have
	\begin{equation}
	\begin{array}{cl}
	r_{k+1} &\leq \left(1-\tfrac{\mu_{\cal D}^2\lambda}{L}\right)^{k+1}r_0 + \tfrac{Lt^2}{8}\sum\limits_{i=0}^{k}\left(1-\tfrac{\mu_{\cal D}^2\lambda}{L}\right)^i\\
	& \leq \left(1-\tfrac{\mu_{\cal D}^2\lambda}{L}\right)^{k+1}r_0 +\tfrac{L^2t^2}{8\mu_{\cal D}^2\lambda}.
	\end{array}
	\end{equation}
	Hence if $t\le \tfrac{2\mu_{\cD}\sqrt{\lambda\e}}{L}$ and $K$ satisfies \eqref{eq:num_it} we get $r_K \leq \varepsilon$.
\end{proof}




\section{Parallel Stochastic Three Points Method}
\label{sec:pa}

Consider the parallel version of \texttt{STP} proposed in Algorithm~\ref{alg:rds_parallel}.

\begin{algorithm}
\caption{{\em Parallel Stochastic Three Points} (\texttt{PSTP})}
\vspace{-4ex}
\label{alg:rds_parallel}
\begin{rm}
\begin{description}
\item[]
\vspace{3ex}
\item[Initialization] \ \\
Choose $x_0\in \R^n$, stepsizes $\alpha_k>0$, parallelism parameter $\tau$,    differentiation stepsize $t_0$.
\vspace{1ex}
\item[For $k=0,1,2,\ldots$] \ \\
\vspace{-2ex}
\begin{enumerate}
\item For $i=1,2,\ldots, \tau$. Generate a random vector $s_{ki} \sim {\cal D}$.
\item Let $s_k = \frac{1}{\tau}\sum\limits_{i=1}^\tau s_{ki}$.
\item Let $x_+ = x_k+\alpha_k s_k$ and $x_- = x_k - \alpha_k s_k$.
\item  $x_{k+1} = \arg \min \{f(x_-), f(x_+),f(x_k)\}$.
\end{enumerate}
\end{description}
\end{rm}
\end{algorithm}

We start our analysis of the complexity in this section  by stating an extra assumption on the probability distribution ${\cal D}$ which is satisfied for all the probability distributions  given in  Lemma~\ref{lem1}.

\begin{ass}\label{ass:P_par} The probability distribution ${\cal D}$ on $\R^n$ satisfies the following properties.
	\begin{enumerate}
		\item If $s_1,s_2\sim {\cal D}$ are independent, then $\Exp[\langle s_1,s_2 \rangle] = 0$.
		\item There exist a constant $\tilde{\mu}_{\cal D} > 0$ and $\tau > 0$ such that if $s_1,s_2,\ldots,s_\tau \sim {\cal D}$ are independent and for all $g\in \R^n$
		$$
		\Exp\left[\left|\left\langle g, \frac{1}{\tau}\sum\limits_{i=1}^\tau s_i  \right\rangle\right|\right] \geq \frac{\tilde{\mu}_{\cal D}}{\sqrt{\tau}}\|g\|_2
		.$$
	\end{enumerate}
\end{ass}
If the first part of the assumption does not hold for distribution ${\cal D}$ we can consider distribution $\bar{\cal D}$ such that $\Exp_{s \sim {\cal D}}[s] = 0$ by adding opposite vector for each vector from ${\cal D}$ and share the probability measure between opposite vectors in equal ratio. The second part of the assumption holds due to Central Limit Theorem for wide range of distributions (this range covers the examples in Lemma~\ref{lem1}) and due to the second part of Lemma~\ref{lem1} we can say that for big enough $\tau$ we have $\tilde{\mu}_{\cal D} \sim \sqrt{\frac{{2}}{{\pi n}}}$ in the case when $\gamma_{\cD} = 1$.

In the next three subsections we will give the adaptation of the main complexity results obtained for \texttt{STP} for \texttt{PSTP}. 
\subsection{Non-convex Case}

The following theorem is the adaptation of Theorem~\ref{thm:nonconvex1}.
\begin{thm}\label{thm:nonconvex_paral} Let Assumptions~\ref{ass:L-smooth},~\ref{ass:P},~\ref{ass:P_sc}~and~\ref{ass:P_par} hold. Choose $\alpha_k=\tfrac{\alpha_0}{\sqrt{k+1}}$, where $\alpha_0>0$. If \begin{equation}\label{eq:isjss8sus_par} K\geq \frac{2\left(\frac{\sqrt{2\tau}(f(x_0)-f_*)}{\alpha_0} + \frac{L\alpha_0}{2\sqrt{\tau}}\right)^2}{\tilde{\mu}_{\cal D}^2 \e^2},\end{equation}
	then  \[\min_{k=0,1,\dots,K} \Exp \left[ \|\nabla f(x_k)\|_2 \right] \leq \e.\]
\end{thm}
\begin{proof}
	By definition of $x_+$ and $x_-$ we have 
\[
			x_+ = x_k + \tfrac{\alpha_k}{\tau}\sum\limits_{i=1}^{\tau}s_{ki}, \qquad 
			x_- = x_k - \tfrac{\alpha_k}{\tau}\sum\limits_{i=1}^{\tau}s_{ki},
\]
	whence
	\begin{equation*}
		\begin{array}{cl}
			f(x_+) &\leq f(x_k) - \alpha_k\langle \nabla f(x_k), \tfrac{1}{\tau}\sum\limits_{i=1}^{\tau}s_{ki}\rangle + \tfrac{L\alpha_k^2}{2\tau^2}\left\|\sum\limits_{i=1}^{\tau}s_{ki}\right\|_2^2,\\
			f(x_-) &\leq f(x_k) + \alpha_k\langle \nabla f(x_k), \tfrac{1}{\tau}\sum\limits_{i=1}^{\tau}s_{ki}\rangle + \tfrac{L\alpha_k^2}{2\tau^2}\left\|\sum\limits_{i=1}^{\tau}s_{ki}\right\|_2^2.
		\end{array}
	\end{equation*}
	Therefore
	\begin{eqnarray*}
		f(x_{k+1}) &\leq  \min\{f(x_+),f(x_-)\} &\leq  f(x_k) - \alpha_k \left|\left\langle\nabla f(x_k), \tfrac{1}{\tau}\sum\limits_{i=1}^\tau s_{ki} \right\rangle\right| + \tfrac{L\alpha_k^2}{2\tau^2}\left\|\sum\limits_{i=1}^{\tau}s_{ki}\right\|_2^2.
	\end{eqnarray*}
	Taking conditional mathematical expectation $\Exp[\;\cdot\mid x_k]$ from the both sides of previous inequality we have
	\begin{equation}\label{eq:par_conv_expected}
		\begin{array}{cl}
			\Exp[f(x_{k+1})\mid x_k] - f_* &\overset{\circledOne}{\leq} f(x_k) - f_* - \tfrac{\alpha_k\tilde{\mu}_{\cal D}}{\sqrt{\tau}}\|\nabla f(x_k)\|_2 + \tfrac{L\alpha_k^2}{2\tau}
		\end{array}
	\end{equation}
	where $\circledOne$ is due to the first part of Assumption~\ref{ass:P_par} and Assumption~\ref{ass:P}:
	\begin{eqnarray*}
		\Exp\left[\left\|\sum\limits_{i=1}^{\tau}s_{ki}\right\|_2^2\right] = \sum\limits_{i=1}^\tau \Exp[\|s_{ki}\|_2^2] + \sum\limits_{i \neq j=1}^{\tau} \Exp[\langle s_{ki},s_{kj} \rangle] = \tau.
	\end{eqnarray*}
	Taking full expectation from the both sides of the inequality \eqref{eq:par_conv_expected} and rearranging the terms we obtain
	\begin{equation}\label{eq:s9jd7d76d_par} g_k \leq \tfrac{\sqrt{\tau}}{\tilde{\mu}_{\cal D}} \left( \tfrac{\theta_k-\theta_{k+1}}{\alpha_k} + \tfrac{L}{2\tau} \alpha_k \right) = \tfrac{\sqrt{\tau}}{\tilde{\mu}_{\cal D}} \left( \tfrac{(\theta_k-\theta_{k+1})\sqrt{k+1}}{\alpha_0} + \tfrac{L\alpha_0}{2\tau \sqrt{k+1}}\right),\end{equation}
	where $g_k = \|\nabla f(x_k)\|_2$. 
	We know from \eqref{eq:monotonicity} and the assumption that $f$ is bounded below that $f_*\leq \theta_{k+1}\leq \theta_k\leq f(x_0)$ for all $k$. Letting $l=\lfloor K/2 \rfloor$, this implies that \[\sum_{j=l}^{2l} \theta_{j}-\theta_{j+1}  = \theta_l - \theta_{2l+1} \leq  f(x_0)-f_*\eqdef C,\]
	from which we conclude that there must exist $j\in \{l,\dots,2l\}$ such that $\theta_j-\theta_{j+1}\leq C/(l+1)$. This implies that
	\begin{eqnarray*}
		g_j &\overset{\eqref{eq:s9jd7d76d_par}}{\leq} & \tfrac{\sqrt{\tau}}{\tilde{\mu}_{\cal D}} \left( \tfrac{(\theta_j-\theta_{j+1})\sqrt{j+1}}{\alpha_0} + \tfrac{L\alpha_0}{2\tau \sqrt{j+1}}\right) \leq  \tfrac{\sqrt{\tau}}{\tilde{\mu}_{\cal D}} \left( \tfrac{C\sqrt{j+1}}{\alpha_0 (l+1)} + \tfrac{L\alpha_0}{2\tau \sqrt{j+1}}\right)\\
		&\leq & \tfrac{\sqrt{\tau}}{\tilde{\mu}_{\cal D}} \left( \tfrac{C\sqrt{2l+1}}{\alpha_0 (l+1)} + \tfrac{L\alpha_0}{2\tau \sqrt{l+1}}\right) \leq \tfrac{\sqrt{\tau}}{\tilde{\mu}_{\cal D} \sqrt{l+1}} \left( \tfrac{\sqrt{2}C}{\alpha_0 } + \tfrac{L\alpha_0}{2\tau}\right)\\
		&\leq & \tfrac{1}{\tilde{\mu}_{\cal D} \sqrt{K/2}} \left( \tfrac{\sqrt{2\tau}C}{\alpha_0 } + \tfrac{L\alpha_0}{2\sqrt{\tau}}\right) \overset{\eqref{eq:isjss8sus_par}}{\leq}  \e.
	\end{eqnarray*}
\end{proof} 

Note that $\alpha_0 = \frac{\sqrt{2 \sqrt{2}\tau}\sqrt{f(x_0)-f_*}}{\sqrt{L}}$ gives the optimal rate which does not depend on $\tau$ and coincides with the rate for spherical setup in the {\tt STP} method. It means that for big enough $\tau$ 
the previous theorem gives a complexity guarantee of the form $O(\frac{n}{\e^2})$.

%

\subsection{Convex Case}
In this subsection we state the main complexity result when $f$ is convex.
\begin{thm}
\label{thm:convex_paral}
 Let Assumptions~\ref{ass:L-smooth},~\ref{ass:levelset} (with $\|\cdot\|_{\cal D} = \|\cdot\|_2$),~\ref{ass:P_sc} and~\ref{ass:P_par} be satisfied. Let $\alpha_k = \alpha_0\left(f(x_k)-f_*\right)$, where $0 < \alpha_0 \leq \frac{2\tau\tilde{\mu}_{\cal D}}{R_0L}$. Define $a  = \frac{\tilde{\mu}_{\cal D}\alpha_0}{\sqrt{\tau}R_0} - \frac{L\alpha_0^2}{2\tau}.$ If $k \ge k(\e) \eqdef \frac{1}{a}\left(\frac{1}{\e} - \frac{1}{r_0}\right),$ then $ \Exp \left[f(x_k) - f(x_*) \right]  \le \e$. 
\end{thm}
\begin{proof}
We have
\begin{equation}\label{eq:par_conv_expe}
		\begin{array}{cl}
			\Exp[f(x_{k+1})\mid x_k] - f_* &\overset{\circledOne}{\leq} f(x_k) - f_* - \tfrac{\alpha_k\tilde{\mu}_{\cal D}}{\sqrt{\tau}}\|\nabla f(x_k)\|_2 + \tfrac{L\alpha_k^2}{2\tau}\\
			&\overset{\circledTwo}{\leq} f(x_k) - f_* - \tfrac{\tilde{\mu}_{\cal D}\alpha_k}{R_0\sqrt{\tau}}(f(x_k) - f_*) + \tfrac{L\alpha_k^2}{2\tau}
		\end{array}
	\end{equation}
	where $\circledOne$ follows from (\ref{eq:par_conv_expected}), and
 $\circledTwo$ follows from Assumption~\ref{ass:levelset}. Using our choice of $\alpha_k = \alpha_0(f(x_k) - f_*)$  and taking full mathematical expectation from the both sides of \eqref{eq:par_conv_expected} we obtain
	\begin{equation*}
		r_{k+1} \leq r_k - \left(\tfrac{\tilde{\mu}_{\cal D}\alpha_0}{\sqrt{\tau}R_0} - \tfrac{L\alpha_0^2}{2\tau}\right)r_k^2 = r_k - ar_k^2.
	\end{equation*}
	Therefore,
$
	\tfrac1{r_{k+1}} 
	-\tfrac1{r_k}
	=
	\tfrac{r_k - r_{k+1}}{r_k r_{k+1} }
	\geq 
	\tfrac{r_k - r_{k+1}}{r_k^2 } 
	\geq a.
$ 
	From this
	we have
$
	\tfrac1{r_k} \geq  \tfrac1{r_0} + ka
$
	and hence
	$
	r_k \le \tfrac1{\tfrac1{r_0} + ka}.
	$
	Finally, if $k \ge \tfrac{1}{a}\left(\tfrac{1}{\e} - \tfrac{1}{r_0}\right)$, then 
	$r_k  \le  \tfrac1{\tfrac1{r_0} + ka} \le  \e$.
\end{proof}
Note that $\alpha_0 = \frac{\sqrt{\tau}\tilde{\mu}_{\cal D}}{R_0L}$ maximizes the value $a$. The optimal value of $a$ is $ \frac{\tilde{\mu}_{\cal D}^2}{2R_0L^2}$, which is proportional to $\frac{1}{\pi n R_0L^2}$ due to the second part of Lemma~\ref{lem1}.
It means that for big enough $\tau$  the above theorem gives an iteration complexity guarantee of the form $O(\frac{n}{\e})$. 

\subsection{Strongly Convex Case}

In this subsection we state the main complexity result when $f$ is strongly convex. The following theorem is an adaptation of Theorem~\ref{thm:stronglyconvex1}.
\begin{thm}
\label{thm:stronglyconvexparal}
	Let Assumptions~\ref{ass:P},~\ref{ass:L-smooth},~\ref{ass:P_sc},~\ref{ass:strongconvex} and~\ref{ass:P_par} be satisfied.
	Let $\alpha_k = \frac{\theta_k \tilde{\mu}_{\cal D}\sqrt{\tau}}{L}\sqrt{2\lambda (f(x_k) - f(x_*) )}$, 
	for some $\theta_k \in (0,2)$ such that $\theta\eqdef \inf_k 2\theta_k-\theta_k^2>0$.
	If    \begin{equation}\label{eq:isjsssus_par}K\geq \frac{L}{\lambda \tilde{\mu}_{\cal D}^2  \theta} \log \left(\frac{f(x_0)-f(x_*)}{\e}\right),\end{equation}
	then $\Exp \left[f(x_K) - f(x_*) \right]  \le \e$.
\end{thm}
\begin{proof}
	By injecting $\alpha_k$ into the first inequality of \eqref{eq:par_conv_expected} we get
	\begin{eqnarray*}
		\Exp[f(x_{k+1}) | x_k] -f(x_*)  &\le& f(x_k) - f(x_*)  -\tfrac{\tilde{\mu}_{\cal D}^2\theta_k \sqrt{2\lambda( f(x_k) - f(x_*) )}\|\nabla f(x_k) \|_2}{L}\\
		&& \qquad + \tfrac{\tilde{\mu}_{\cal D}^2\theta_k^2 \lambda (f(x_k) - f(x_*) )}{L}.
	\end{eqnarray*}
	From the strong convexity property of $f$ we have $\|\nabla f(x_k) \|_2^2 \ge 2\lambda (f(x_k) - f(x_*) ).$ Therefore
	\begin{eqnarray*}
		\Exp[f(x_{k+1}) | x_k] -f(x_*)  &\le& f(x_k) - f(x_*)  -\tfrac{2\tilde{\mu}_{\cal D}^2\theta_k \lambda( f(x_k) - f(x_*) )}{L} + \tfrac{\tilde{\mu}_{\cal D}^2\theta_k^2 \lambda (f(x_k) - f(x_*) )}{L}\\
		&\le & f(x_k) - f(x_*)  -\tfrac{\tilde{\mu}_{\cal D}^2 \lambda( f(x_k) - f(x_*) )}{L} (2\theta_k - \theta_k^2)\\
		&\le & f(x_k) - f(x_*)  -\tfrac{\tilde{\mu}_{\cal D}^2\theta \lambda( f(x_k) - f(x_*) )}{L},
	\end{eqnarray*}
by using the definition of  $\theta$.
	Let $r_k = \Exp \left[f(x_k)\right] - f(x_*) $. By taking the expectation of the last inequality we get
$
		r_{k+1}   \le \left(1 - \tfrac{\tilde{\mu}_{\cal D}^2\theta\lambda}{L}\right) r_k
$,   
	hence
	\begin{eqnarray*}
		r_{k}   \le \left(1 - \tfrac{\tilde{\mu}_{\cal D}^2\theta\lambda}{L}\right)^k r_0. 
	\end{eqnarray*}
	Therefore, if $K$ satisfies \eqref{eq:isjsssus_par}, we get $r_K \le \e$.
\end{proof}

For big enough $\tau$  the above theorem gives an iteration complexity guarantee of the form $O(n\log\left(\tfrac{1}{\e}\right))$. 

\section{Numerical Results} 
\label{sec:exper}
In this section, we report the results of some preliminary experiments performed in order to assess the efficiency and
the robustness of the proposed algorithms compared to the coordinate search method (this method will be called \texttt{\texttt{DDS}}) and the algorithm proposed in \cite{Nesterov_Spokoiny_2017}. 
In the latter approach, 
 at each iteration $k$,  a random vector $s_k$ following the uniform distribution on the unit sphere is generated,  
then the next iterate is computed as follows  
\begin{equation}
\label{eq:nesupdate}
x_{k+1} = x_k - \alpha_k \frac{f(x_k + \mu_k s_k) - f(x_k)}{ \mu_k} s_k, 
\end{equation}
where $\mu_k\in (0,1)$ is the finite differences parameter, and $\alpha_k$ is the stepsize.
This method generates a trial step similar to one of the trial steps in our method ($x_- = x_k - \alpha_k s_k$)  when the probability distribution ${\cal D}$ is chosen to be the uniform distribution on the unit sphere up to the multiplication 
of the step by $\frac{f(x_k + \mu_k s_k) - f(x_k)}{ \mu_k}$.  
 This  method will be called \texttt{RGF} (Random Gradient free method).
 All the results presented here are averaged over 10 runs of the algorithms.  
We did all our experiments using Matlab. 

To compare the performance of the algorithms we use performance profiles proposed by Dolan and Mor\'e~\cite{Dolan_2002} over a variety of problems.  Given a set of problems $\mathcal{P}$ (of cardinality $|\mathcal{P}|$)
and a set of algorithms (solvers) $\mathcal{S}$, the performance profile
$\rho_s(\tau)$ of an algorithm~$s$ is defined as the fraction of problems
where the performance ratio $r_{p,s}$ is at most $\tau$
\begin{eqnarray*}
 \rho_s(\tau) \; = \; \frac{1}{|\mathcal{P}|} \mbox{size} \{ p \in \mathcal{P}: r_{p,s} \leq \tau \}.
\end{eqnarray*}
The performance ratio $r_{p,s}$ is in turn defined by
\[
r_{p,s} \; = \; \frac{t_{p,s} }{\min\{t_{p,s}: s \in \mathcal{S}\}},
\]
where $t_{p,s} > 0$ measures the performance of the algorithm~$s$ when solving problem~$p$, seen here as the number of function evaluation.
Better performance of the algorithm~$s$,
relatively to the other algorithms on the set of problems,
is indicated by higher values of $\rho_s(\tau)$.
In particular, efficiency is measured by $\rho_s(1)$ (the fraction of problems for which algorithm~$s$ performs the best) and robustness is measured by $\rho_s(\tau)$ for $\tau$ sufficiently large
(the fraction of problems solved by~$s$). Following what is suggested in~\cite{Dolan_2002} for a better visualization,
we will plot the performance profiles in a $\log_2$-scale
(for which $\tau=1$ will correspond to $\tau=0$).

The distribution ${\cal D}$ used here for our random direction generation is the uniform distribution on the unit sphere. We performed other experiments (not reported here) with  different choices for distributions ${\cal D}$. For instance, the distributions listed in Lemma~\ref{lem1}. 
We found similar performance as those reported here.
The parameters defining the implemented algorithms are set as follows: 
 For \texttt{RGF} we choose $\mu_k = 10^{-4}$, and $\alpha_k = \frac{1}{4(n+4)}$ where $n$ is the problem dimension. For this method the authors proposed to use the stepsize $\alpha_k = \frac{1}{4L(n+4)}$, where $L$ is the Lipschitz constant of the gradient of the objective function. Since for our test problems we do not know this constant, we ran \texttt{RGF}
 method with different values for $L$, for instance $0.1, ~1,~10,$ and $100$. The best performance was found for $L=1$. The stepsize in \texttt{DDS} is initialized by $\alpha_0 =1$, then it is updated dynamically with the iterations by multiplying it by $2$ when the step is successful and dividing it by $2$ otherwise. 


For all algorithms,  %
we counted the number of function evaluations taken to (i) drive the function value
below $f^* + \e \left(f(x_0) - f^* \right)$, where $f^*$  is a local minimal value of the objective function $f$, and $\e$ is a tolerance. In our experiments $\e = 10^{-1}, ~10^{-3}$ and $10^{-5}$, (ii) or the maximum number of iterations  attains $100000$.

\subsection{Non-Convex Case}

In this section, we report the results of comparison of our approach \texttt{STP} for non-convex problems with \texttt{\texttt{DDS}} and \texttt{RGF}. 
We will call our \texttt{STP} method when using the variable stepsize \texttt{STP-vs}, and \texttt{STP-fs} when we use a fix stepsize. For \texttt{STP-vs} we choose $ \alpha_k =\frac{1}{\sqrt{k+1}}$. For 
 \texttt{STP-fs} we choose $\alpha_k = \alpha = 0.1 \e$.
 
We use the Mor\'e/Garbow/Hillstrom 34 test problems \cite{More_1981} which are implemented in Matlab. All the test problems are smooth. The dimension $n$ of the problems changes between $n=2$ to $n=100$, typically  $n=2, 10, 50$ and $100$.
 We use the starting points and the values $f^*$ suggested in~\cite{More_1981} for all the problems. 

\begin{figure}[!ht]
\centering
\subfigure[$\e = 10^{-1}$]{
\includegraphics[scale=0.205]{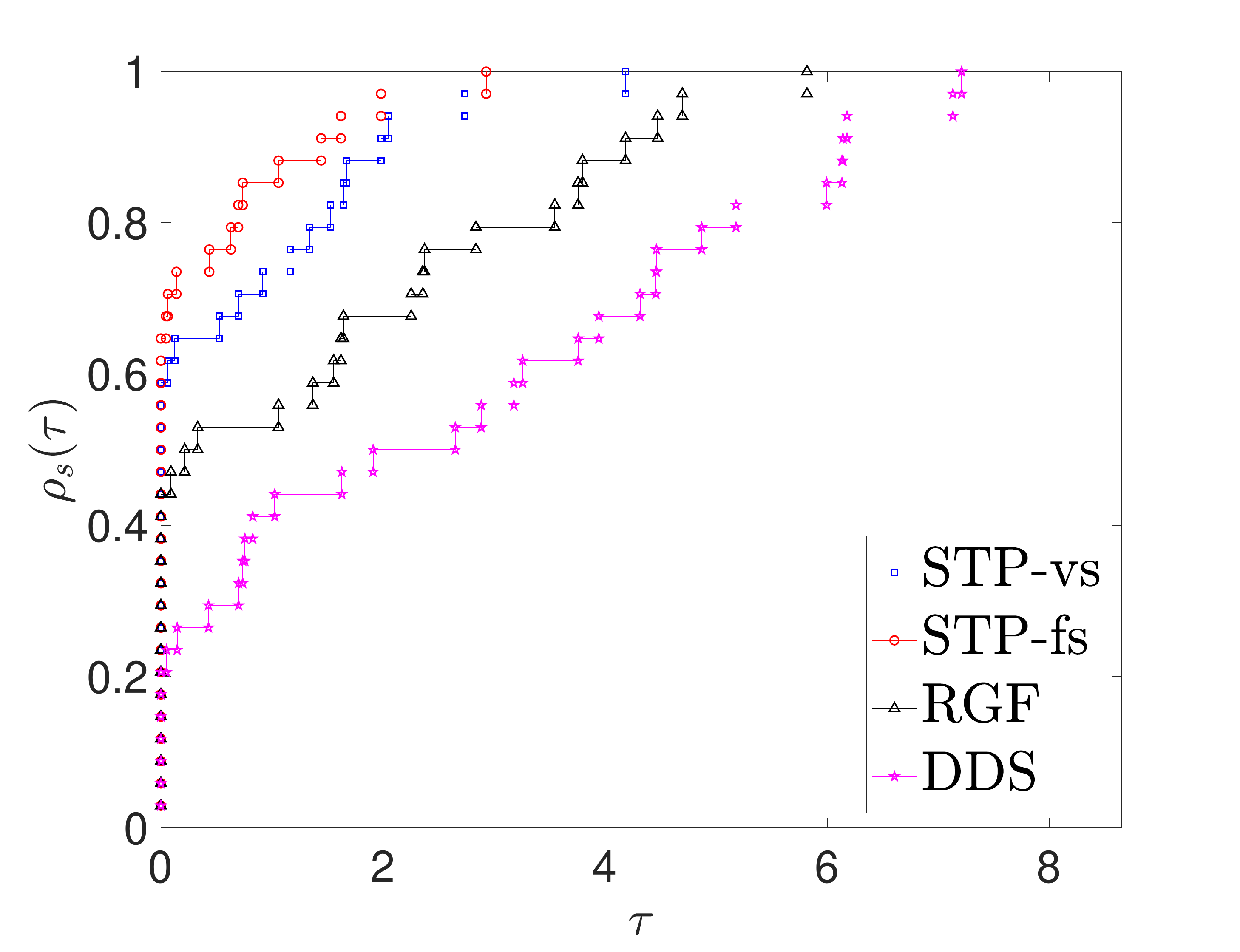}\label{subfig1:bigeps}
}
\subfigure[$\e = 10^{-3}$]{
\includegraphics[scale=0.205]{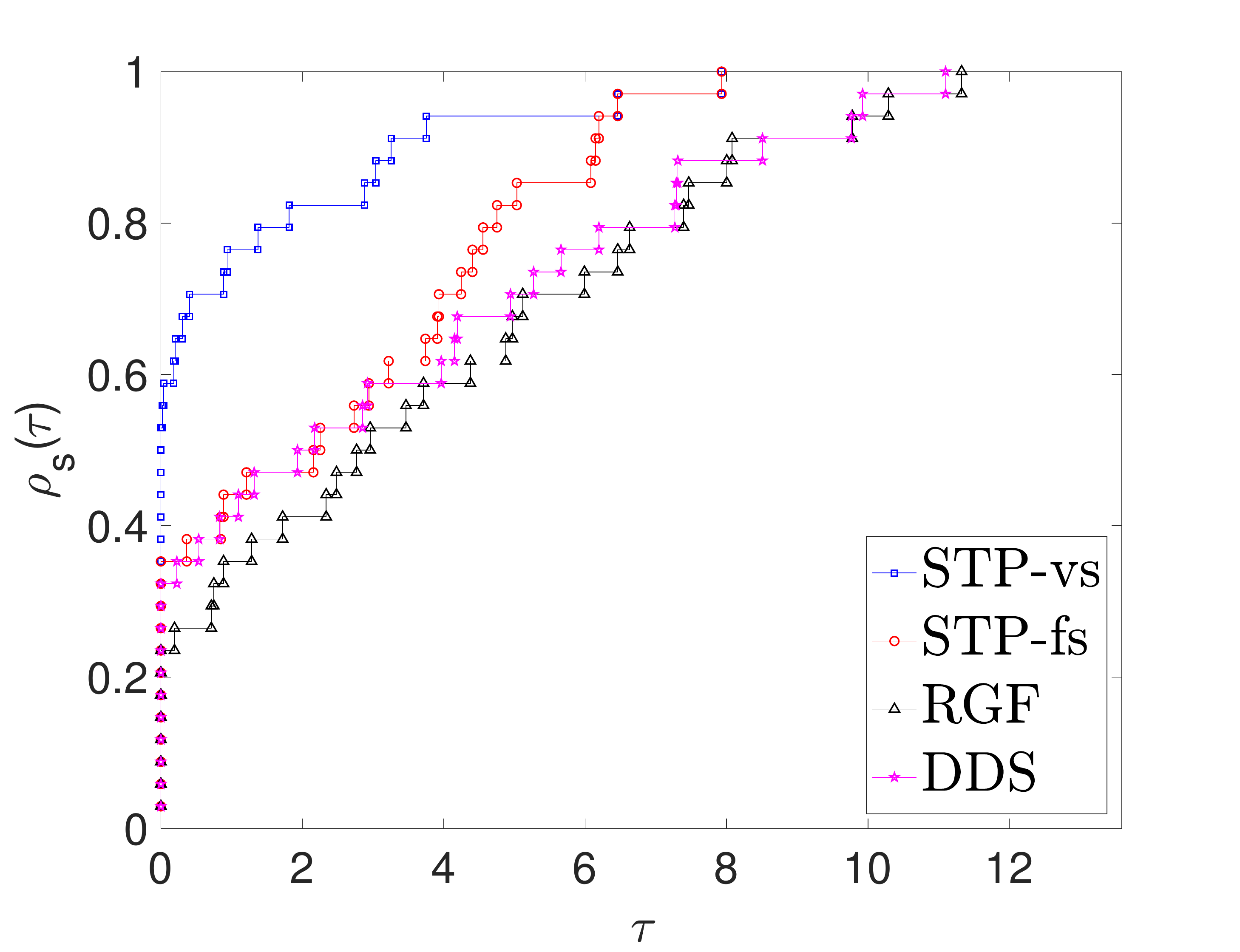} \label{subfig2:medeps}
}
\subfigure[$\e = 10^{-5}$]{
\includegraphics[scale=0.205]{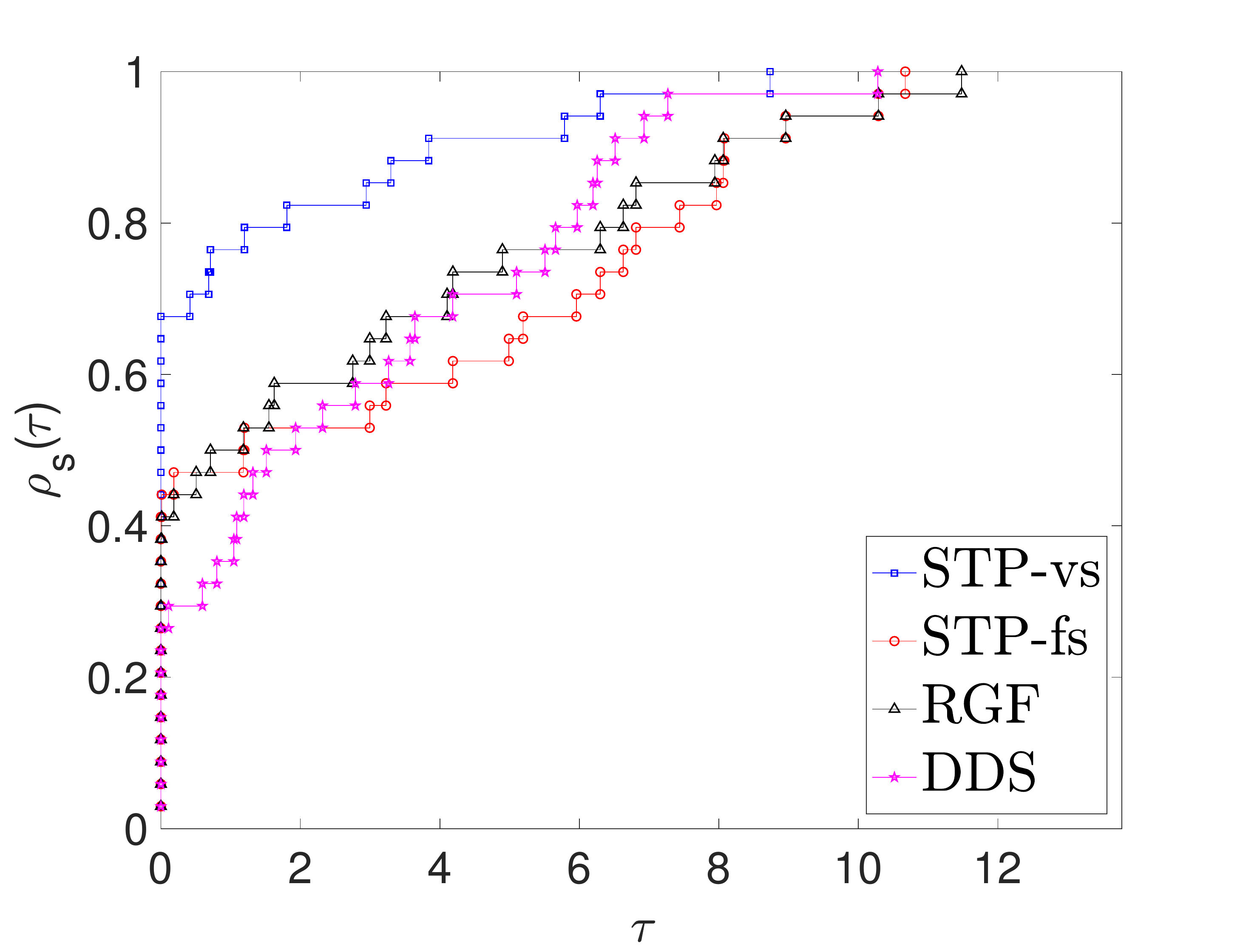} \label{subfig3:smalleps}
}
\caption{Performance profiles on $34$ optimization problems.} \label{fig:ncc}
\end{figure}

Figure~\ref{fig:ncc} depicts the performance profiles of the algorithms. It shows that our approach
 (the methods \texttt{STP-vs} and \texttt{STP-fs}) improves the efficiency  of the \texttt{DDS} and \texttt{RGF} algorithms on the tested problems. In fact, the number of the function evaluations performance profiles show that the use of the random directions leads to a significant improvement on terms of the efficiency (for $\tau=0$, on about $40\%$ of the tested problems our approach performs the best,  and less than $5\%$ for \texttt{RGF} and \texttt{\texttt{DDS}}). From Figures~\ref{subfig1:bigeps} and~\ref{subfig2:medeps}, we see that the use of the random directions leads to a better robustness when a  small precision is targeted (i.e. $\e =10^{-1}$ and $\e =10^{-3}$). 
However, when a big precision ( $\e =10^{-5}$) is targeted  \texttt{DDS} becomes competitive. In fact, as shown in  Figure~\ref{subfig3:smalleps}, \texttt{DDS} is more robust than 
\texttt{RGF} approach and our method using fix stepsize. Our method  \texttt{STP-vs} still more robust than \texttt{DDS}.
 
\subsection{Convex Case}

In this section, we report the results of comparison of two \texttt{STP} methods for convex problems with \texttt{\texttt{DDS}} and \texttt{RGF}. The first \texttt{STP} method is the one using 
the variable stepsize $ \alpha_k =\left| \tfrac{1}{Lt}(f(x_k + t s_k) - f(x_k))\right|$, where $t = 10^{-4}$. 
We will call this method \texttt{STP-vs}. The second \texttt{STP} method is the one using the fix stepsize $\alpha_k = \alpha = 0.1\e$. It will be called 
 \texttt{STP-fs}.

We selected from the Mor\'e/Garbow/Hillstrom  problems those with a unique minimum.  To have a large bed test, we create different instances for problems by varying the problem dimension $n$ when it is possible. Our test bed in this section contains $40$ problems.

\begin{figure}[!ht]
\centering
\subfigure[$\e = 10^{-1}$]{
\includegraphics[scale=0.205]{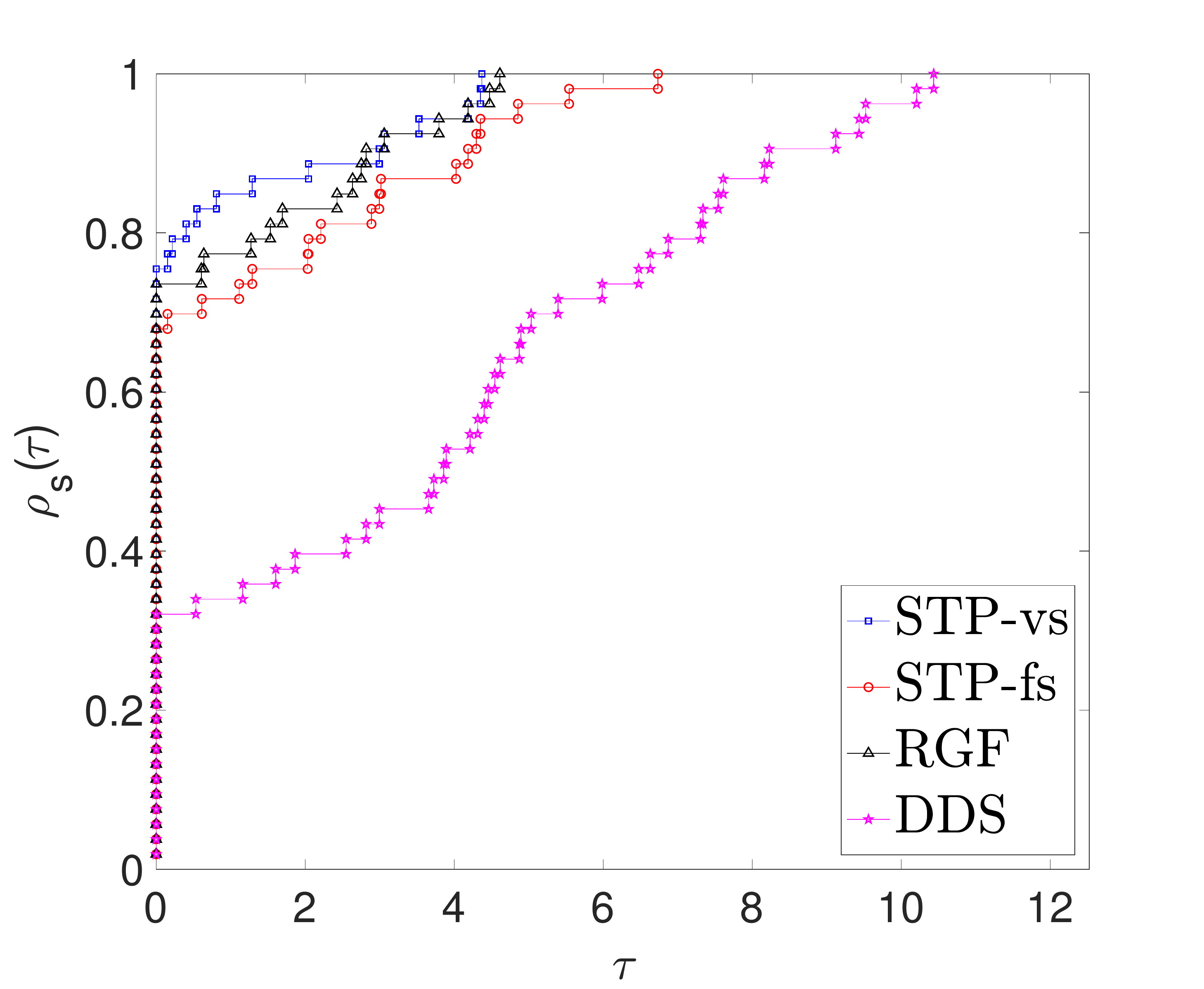}\label{subfig1:bigepsconv}
}
\subfigure[$\e = 10^{-3}$]{
\includegraphics[scale=0.205]{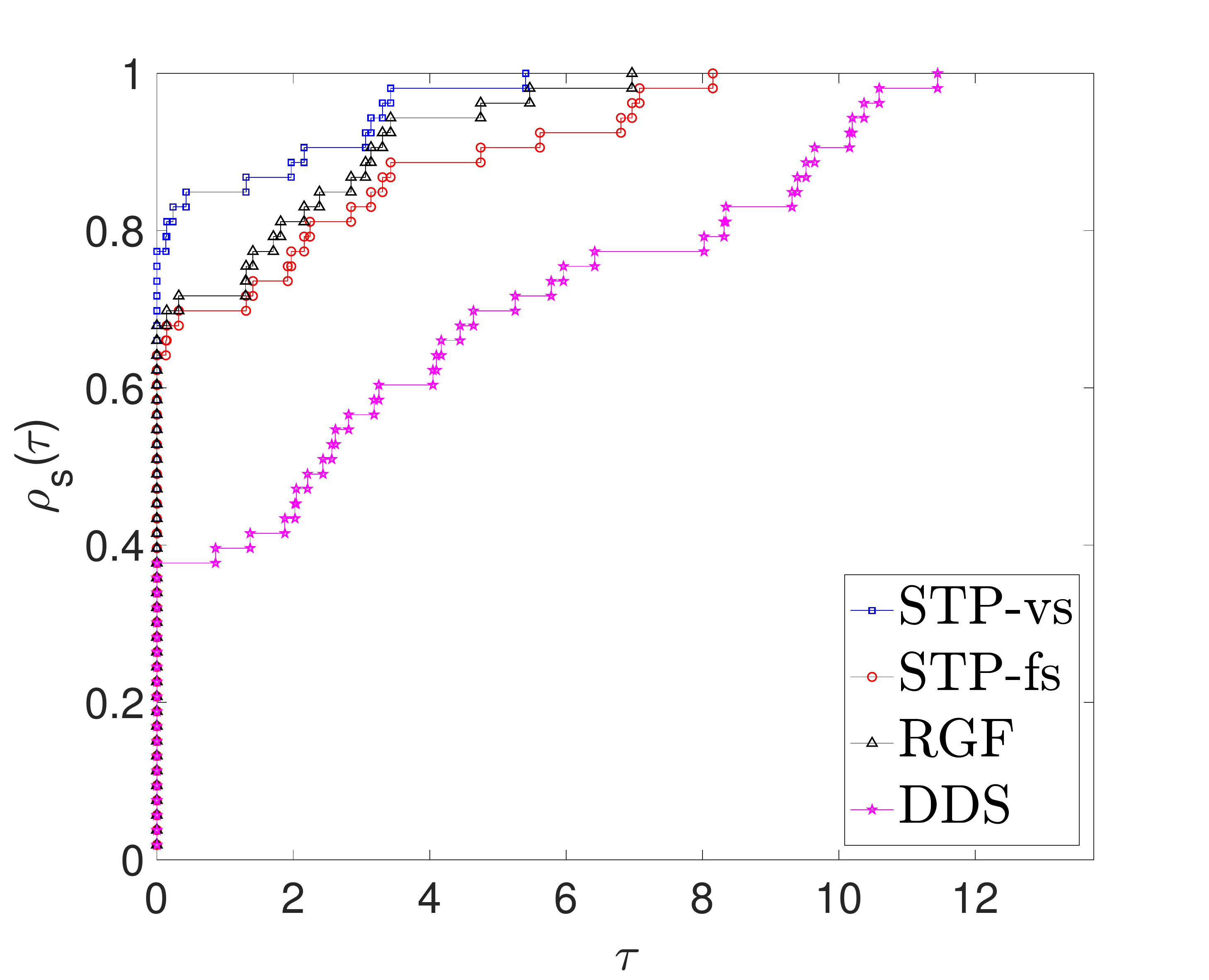} \label{subfig2:medepsconv}
}
\subfigure[$\e = 10^{-5}$]{
\includegraphics[scale=0.205]{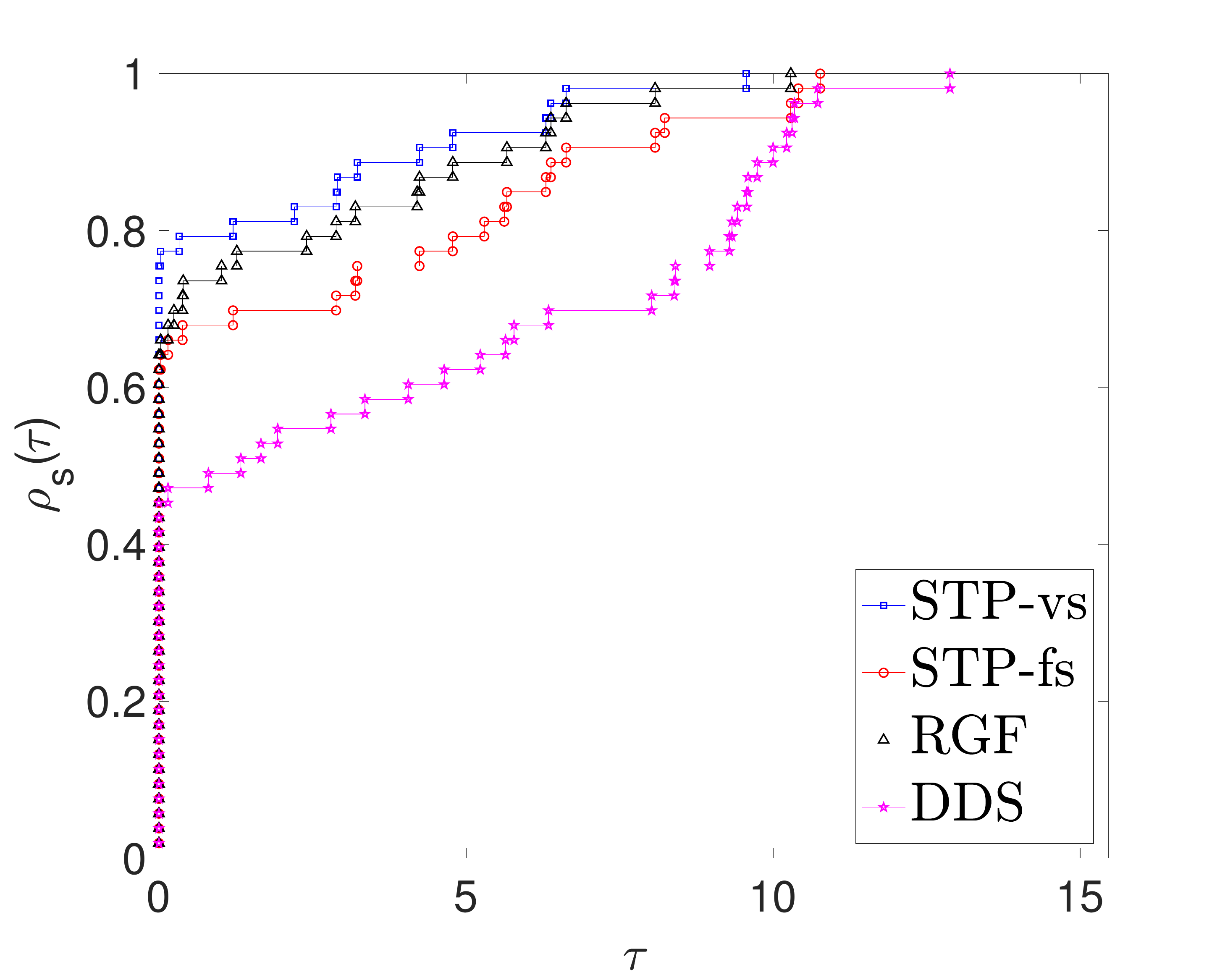} \label{subfig3:smallepsconv}
}
\caption{Performance profiles on $40$ optimization problems.} \label{fig:cc}
\end{figure}

In Figure~\ref{fig:cc}, the performance profiles show that the random based methods (
\texttt{RGF} method and our two methods \texttt{STP-vs} and \texttt{STP-vs})
outperform by far the \texttt{DDS} method. Our method \texttt{STP-vs} gives the best performances for small precision (see Figures~\ref{subfig1:bigepsconv} and~\ref{subfig2:medepsconv}). For big precision ($\e = 1e-5$), it gives almost similar performances as \texttt{RGF} method ((see Figure~\ref{subfig3:smallepsconv}).
Our method \texttt{STP-fs} is outperformed by \texttt{RGF}. 

\subsection{First order methods}

In this section, we report the results of comparison of gradient based methods that our approach cover, using the variable stepsize $ \alpha_k =\frac{1}{\sqrt{k+1}}$ and the fix stepsize  $ \alpha_k =0.1 \e$. In fact to select these stepsizes, we run many experiments with different values and found the best results for the chosen stepsizes.
 We denote with \texttt{ngd-vs}, and \texttt{ngd-fs}  the Normalized Gradient Descent (NGD) methods using the variable stepsize, and the fix stepsize respectively.  
 With similar notation we denote by \texttt{signgd-vs}, and \texttt{signgd-fs} the Signed Gradient Descent (SignGD) methods and by \texttt{nrcd-vs}, and \texttt{nrcd-fs}
  Normalized Randomized Coordinate Descent (NRCD) methods using the  variable stepsize, and the fix stepsize respectively.
  
We use the Mor\'e/Garbow/Hillstrom 34 test problems for which we add 20 problems by creating different instances for problems by varying the problem dimension $n$ when it is possible. Our test bed in this section contains $54$ problems. 

\begin{figure}[!ht]
\centering
\subfigure[$\e = 10^{-1}$]{
\includegraphics[scale=0.205]{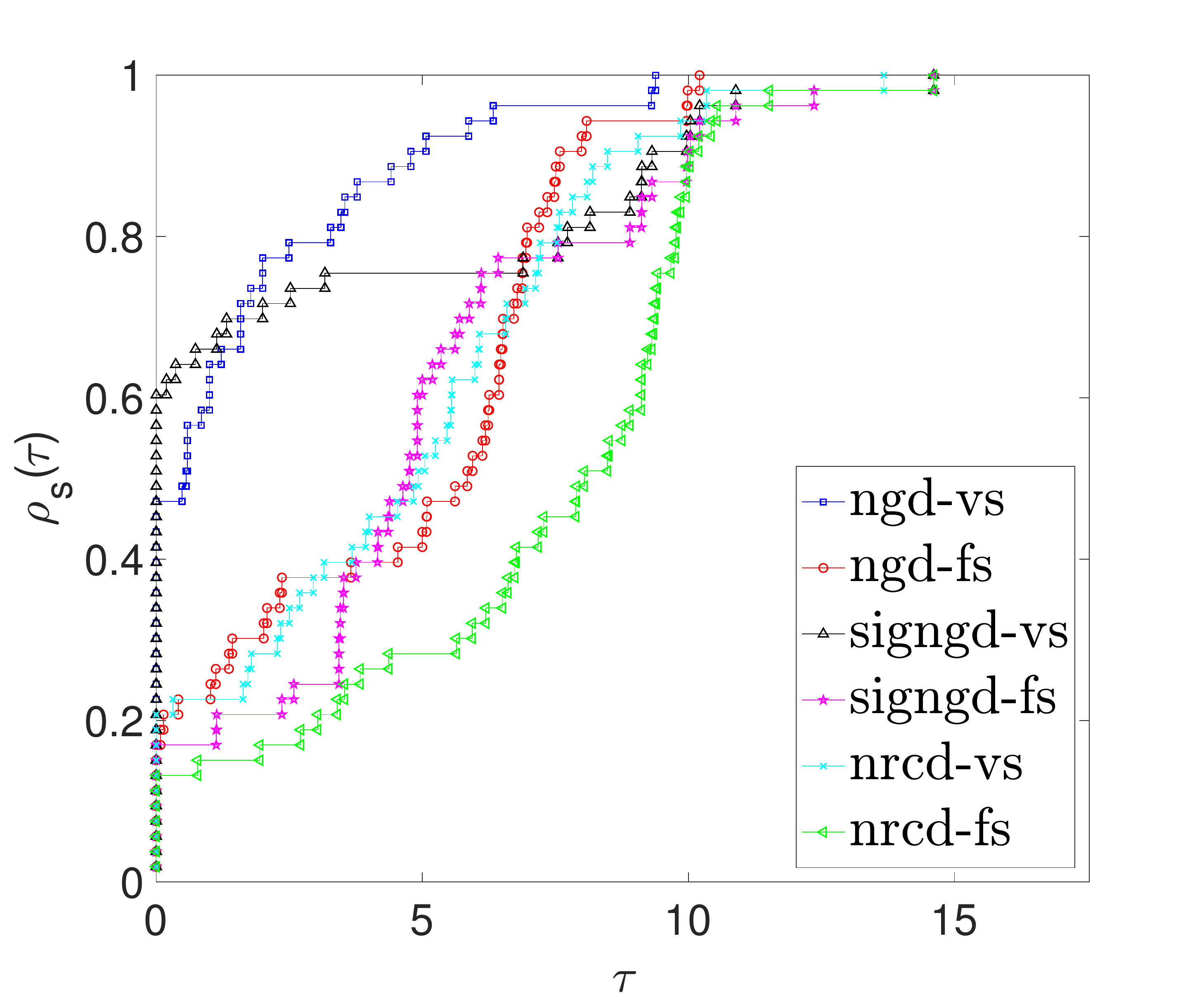}\label{subfig1:bigepsfirstorder}
}
\subfigure[$\e = 10^{-3}$]{
\includegraphics[scale=0.205]{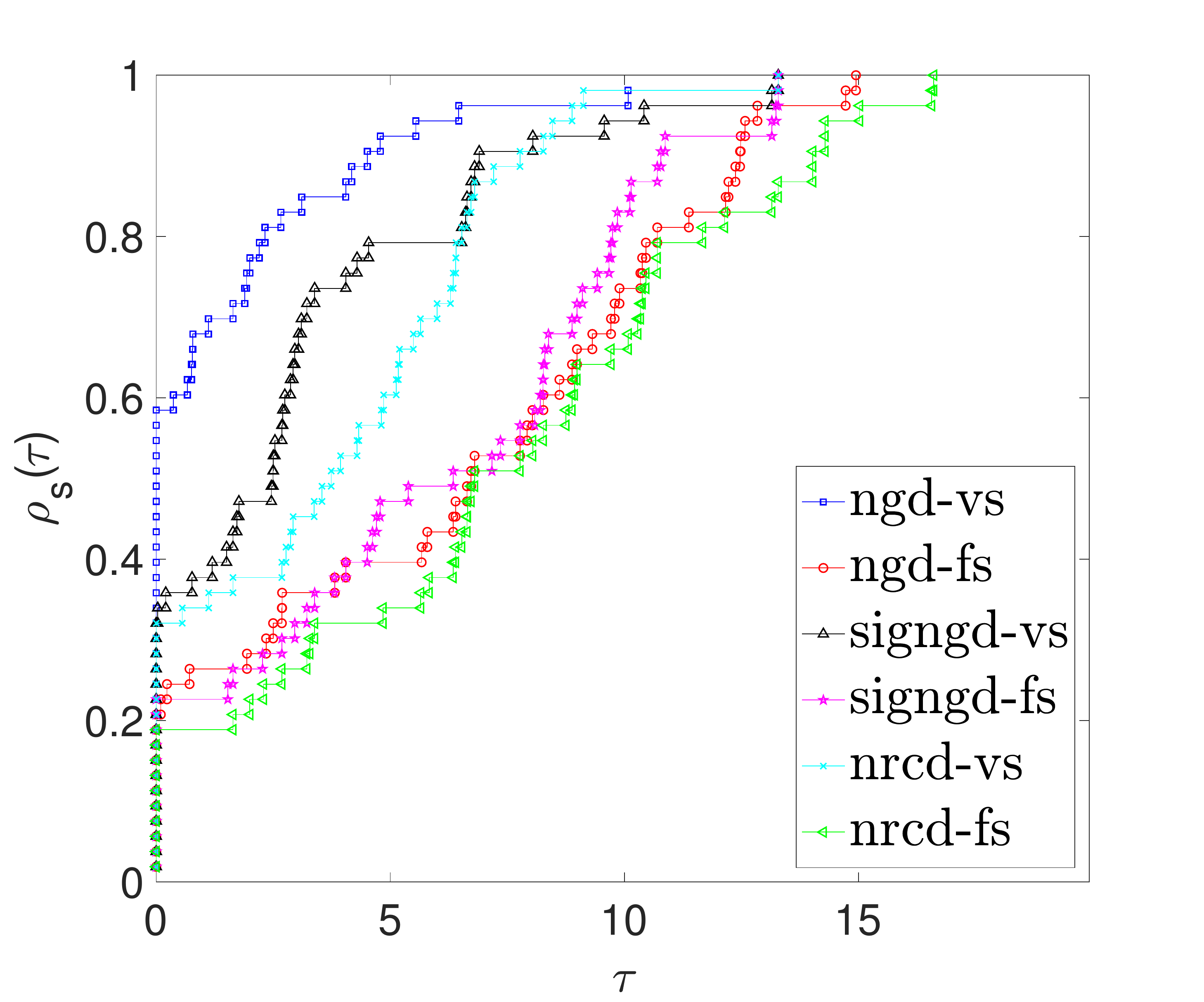} \label{subfig2:medepsfirstorder}
}
\subfigure[$\e = 10^{-5}$]{
\includegraphics[scale=0.205]{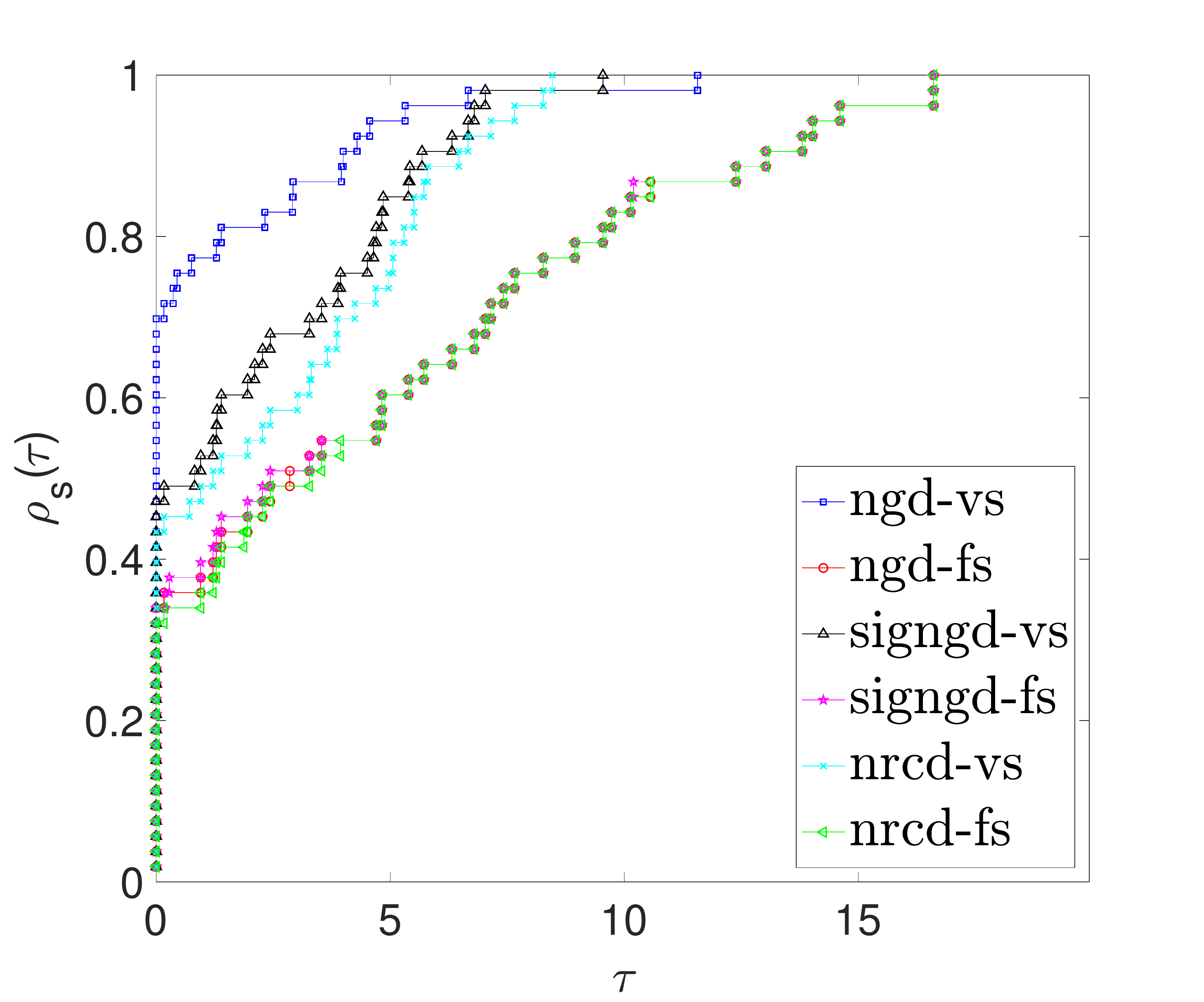} \label{subfig3:smallepsfirstorder}
}
\caption{Performance profiles on $54$ optimization problems.} \label{fig:fo}
\end{figure}

Figure~\ref{fig:fo} depicts the performance profiles of the algorithms. It shows that the use of the variable stepsize gives better performances than the fix stepsize.
As one may expect, the normalized gradient descent method \texttt{ngd-vs} exhibits performances better than the other methods, except for small precision ($\e = 1e-1$), it is less efficient than signed gradient descent method \texttt{signgd-vs}. The  latter method 
 is more efficient and less robust than normalized randomized coordinate descent method \texttt{nrcd-vs}. 
  
\subsection{Experiments for \texttt{PSTP}} 
\begin{figure}[!ht]
	\centering
	\subfigure[$n = 25$]{
		\includegraphics[scale=0.145]{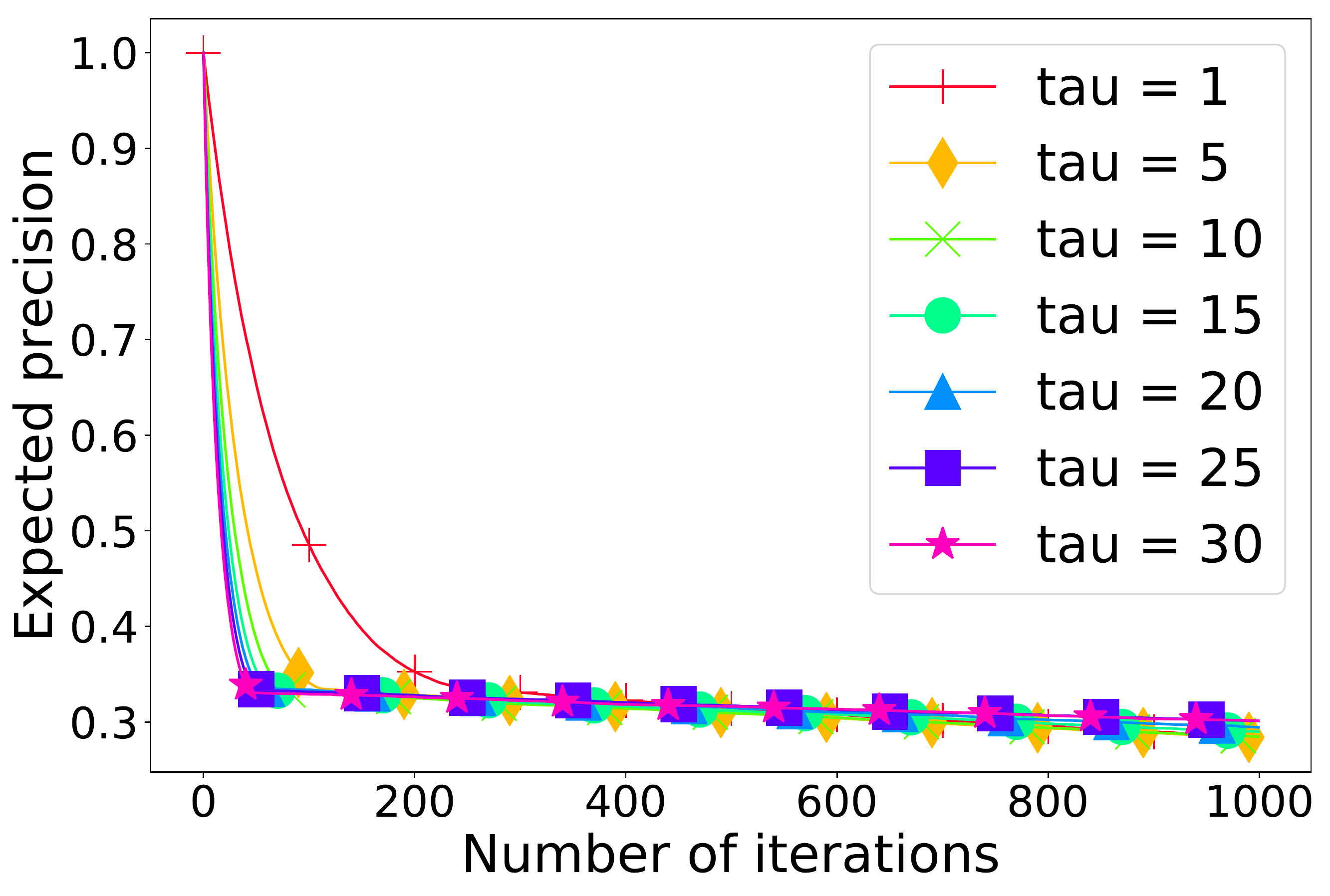}\label{subfig1:n25_par}
	}
	\subfigure[$n = 50$]{
		\includegraphics[scale=0.145]{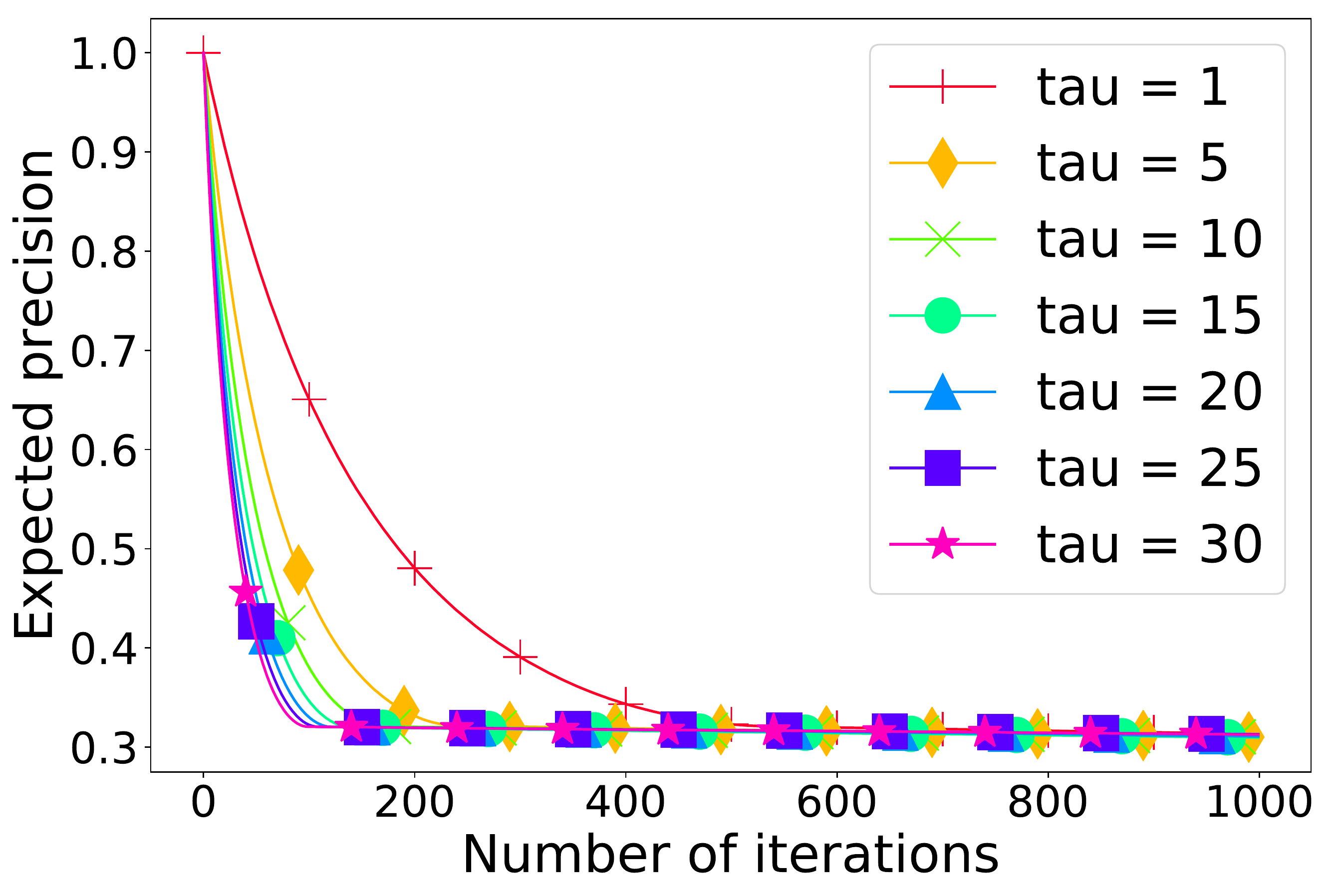}\label{subfig1:n50_par}
	}
	\subfigure[$n = 75$]{
		\includegraphics[scale=0.145]{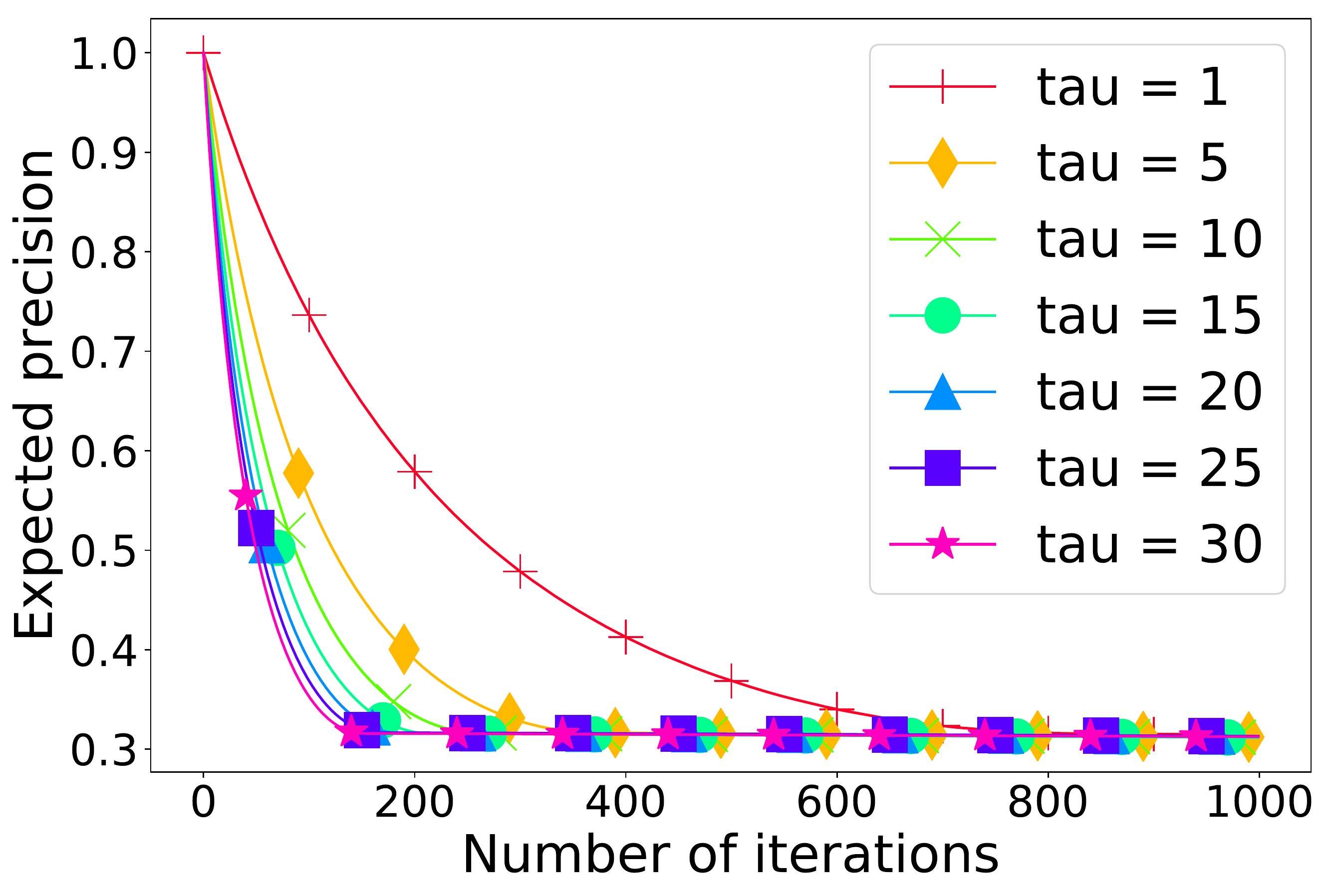}\label{subfig1:n75_par}
	}
	\subfigure[$n = 100$]{
		\includegraphics[scale=0.145]{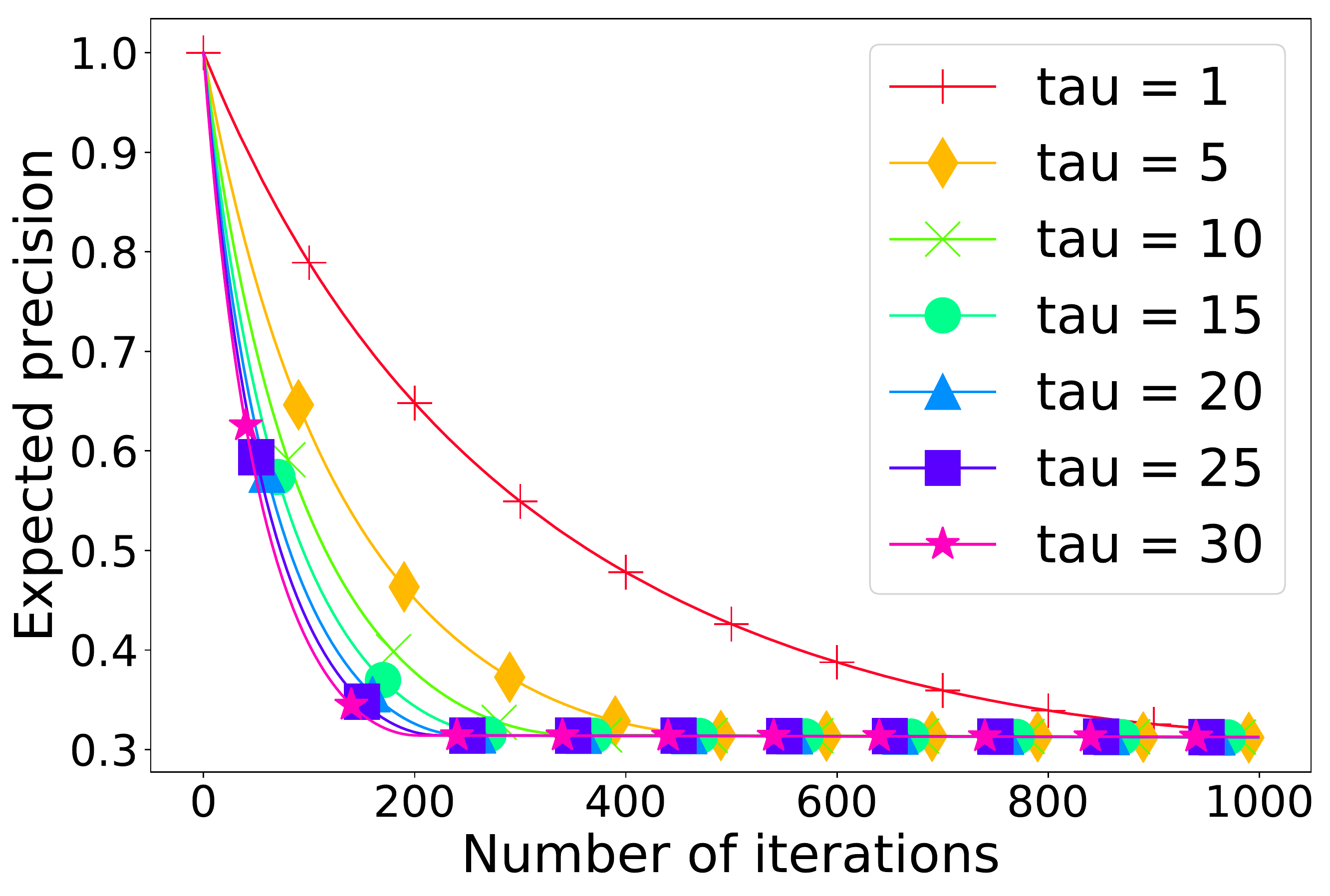}\label{subfig1:n100_par}
	}
	\subfigure[$n = 125$]{
		\includegraphics[scale=0.145]{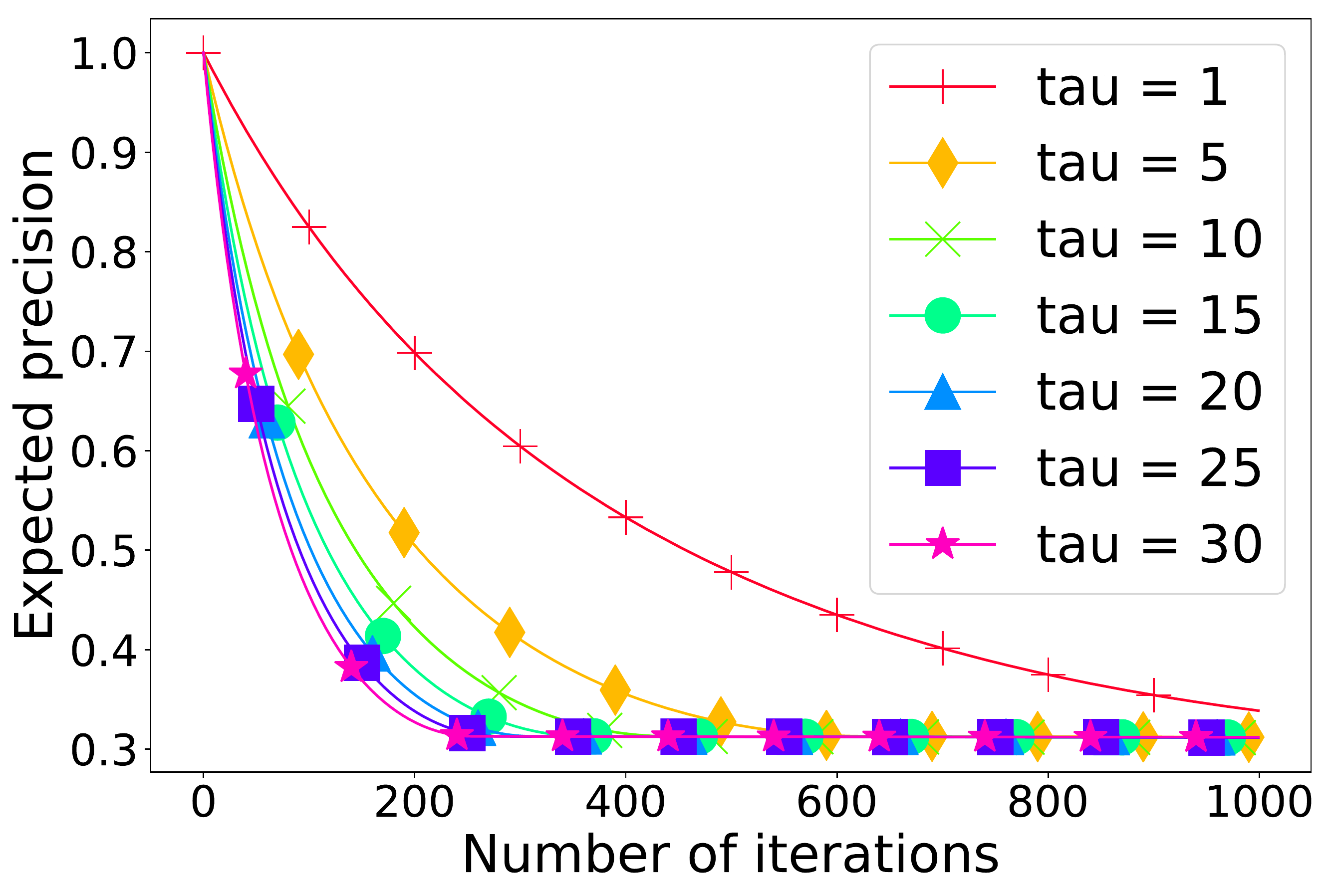}\label{subfig1:n125_par}
	}
	\subfigure[$n = 150$]{
		\includegraphics[scale=0.145]{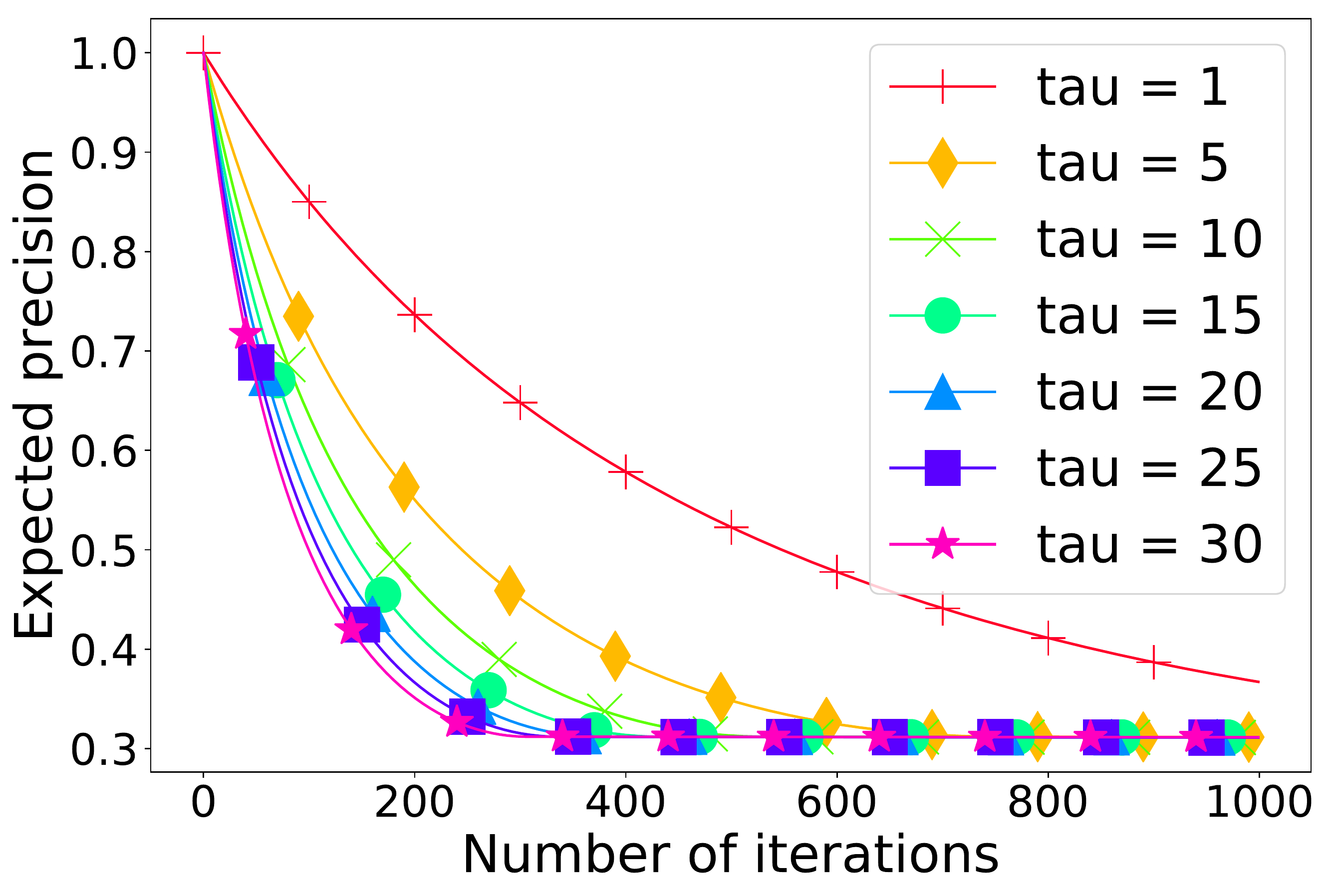}\label{subfig1:n150_par}
	}
	\subfigure[$n = 175$]{
		\includegraphics[scale=0.145]{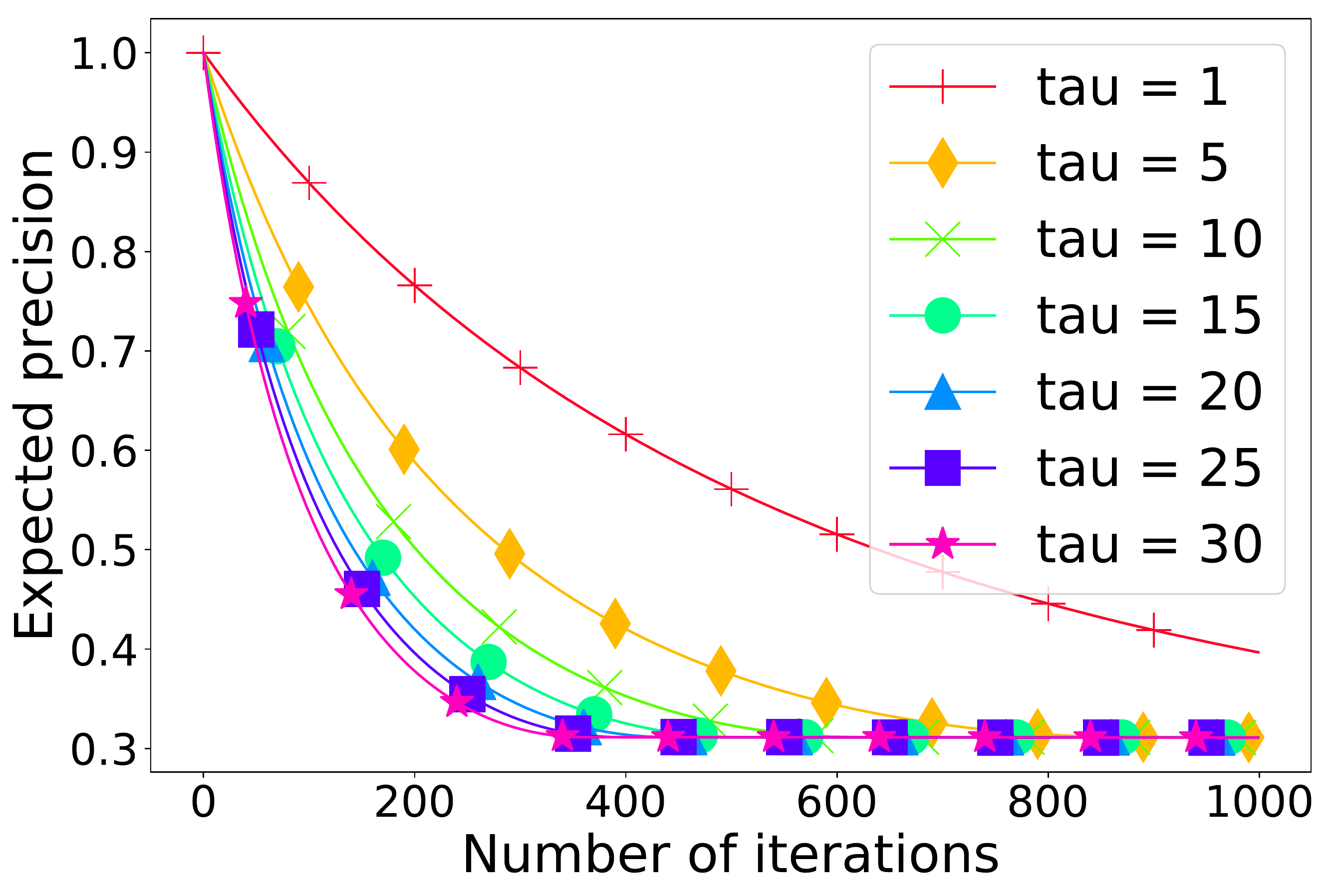}\label{subfig1:n175_par}
	}
	\subfigure[$n = 200$]{
		\includegraphics[scale=0.145]{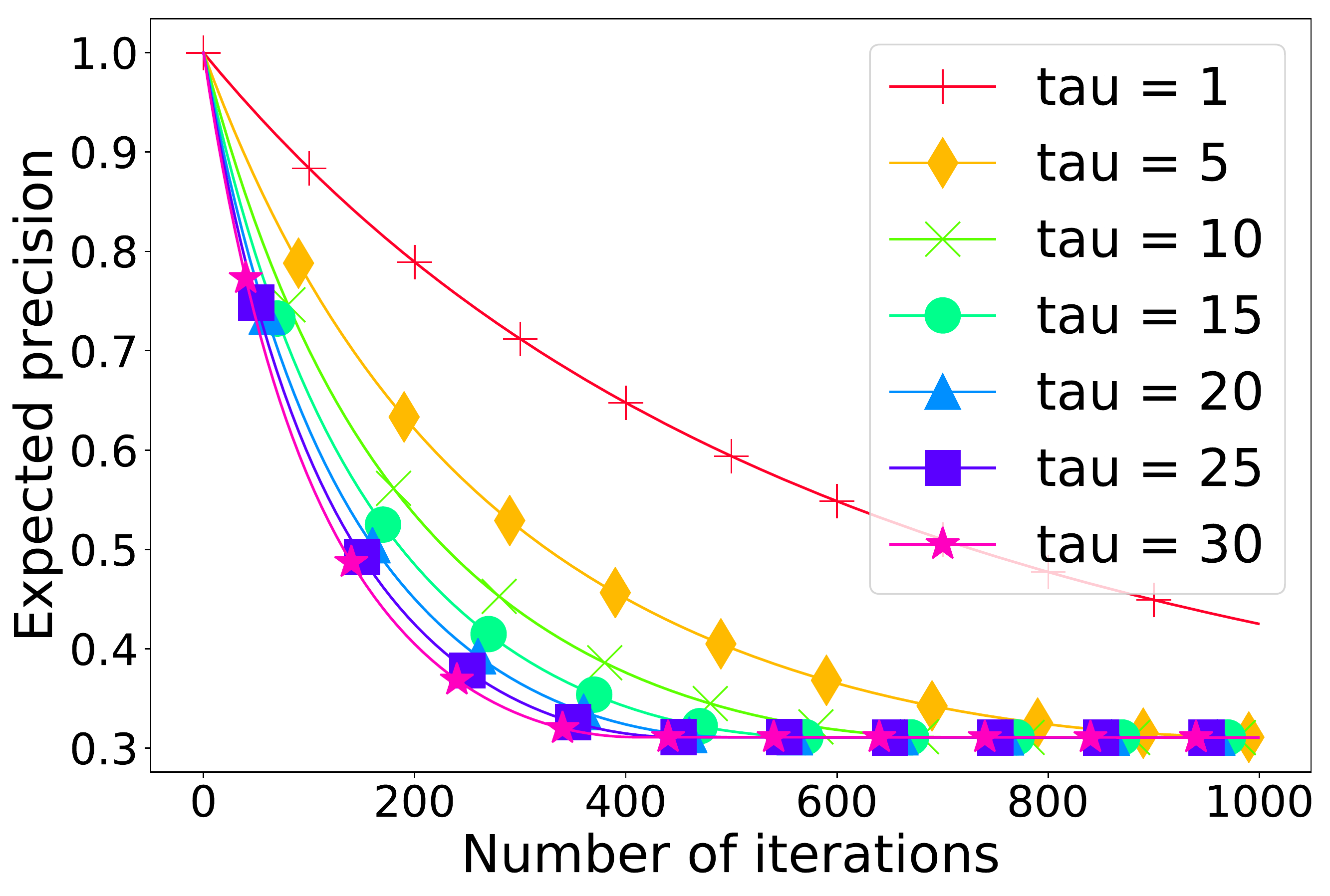}\label{subfig1:n200_par}
	}
	\caption{Trajectories of \texttt{PSTP} for the different $n$.} \label{fig:par}
\end{figure}

We considered the following function
$$
f(x) = \frac{1}{2}x_1^2 + \frac{1}{2}\sum\limits_{i=1}^{n-1}(x_{i+1}-x_i)^2 + \frac{1}{2}x_n^2 - x_1
$$
and run \texttt{PSTP} with different $\tau$ (see Figure~\ref{fig:par}). From the numerical results we see that the rate could be worse for small $\tau$ than for $\tau$. It happens because the Assumption~\ref{ass:P_sc} does not have to be true for small $\tau$ with the parameter $\tilde{\mu}_{\cal D} \sim \sqrt{\frac{{2}}{{\pi n}}}$ as we use in the experiments (recall that this parameter corresponds to the statement of Central Limit Theorem and, therefore, we need $\tau$ to be big enough). When $\tau$ is big enough the Assumption~\ref{ass:P_par} will holds and we will obtain the improvement of the rate. We measure $\frac{f(x_k) - f_*}{f(x_0)-f_*}$ on the $y$-axis and call it ``Expected precision". 
 
\subsection{\texttt{STP} vs \texttt{RGF}} 
\begin{figure}[!ht]
	\centering
	\subfigure[$n = 25$]{
		\includegraphics[scale=0.145]{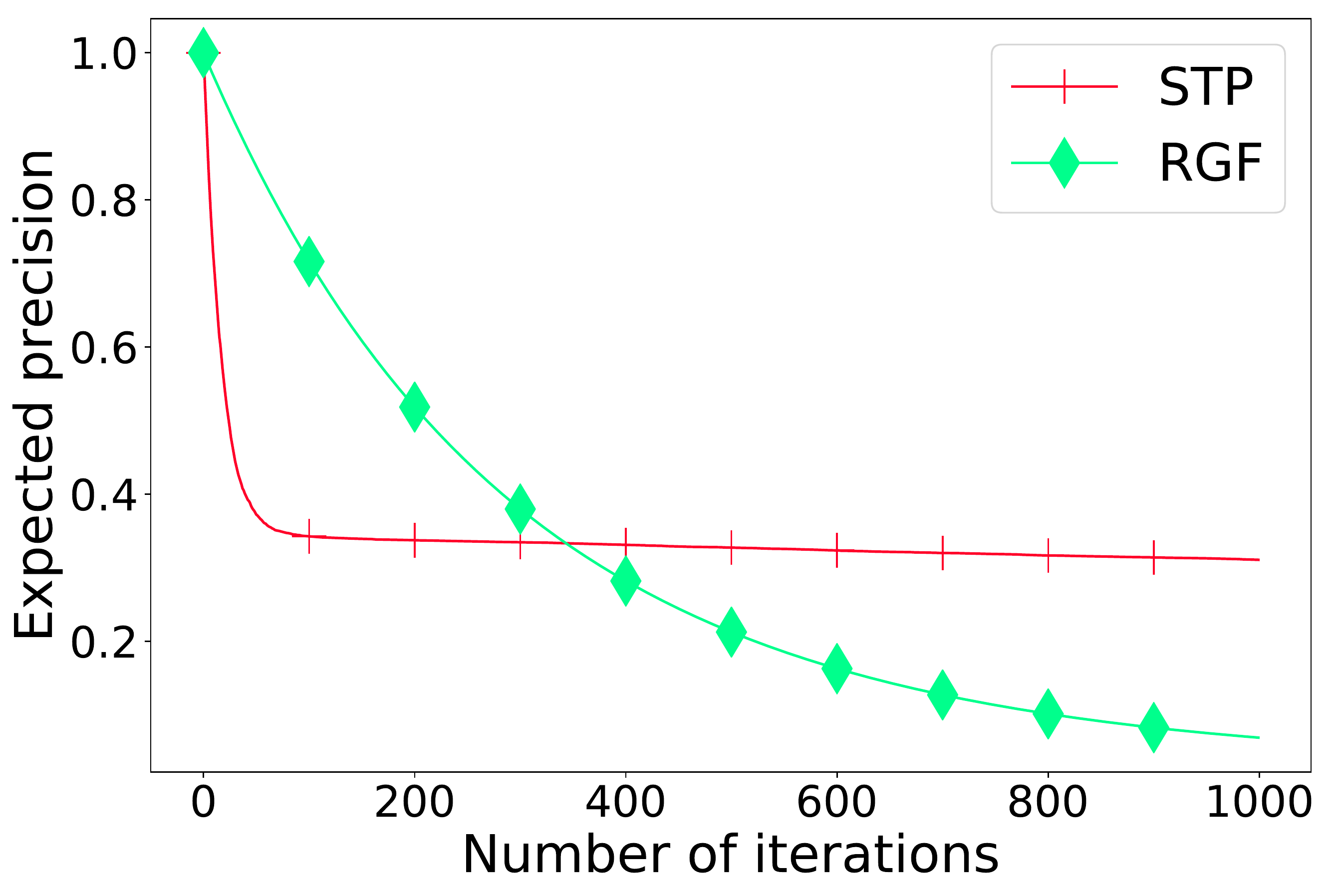}\label{subfig1:n25_rds_rgf}
	}
	\subfigure[$n = 50$]{
		\includegraphics[scale=0.145]{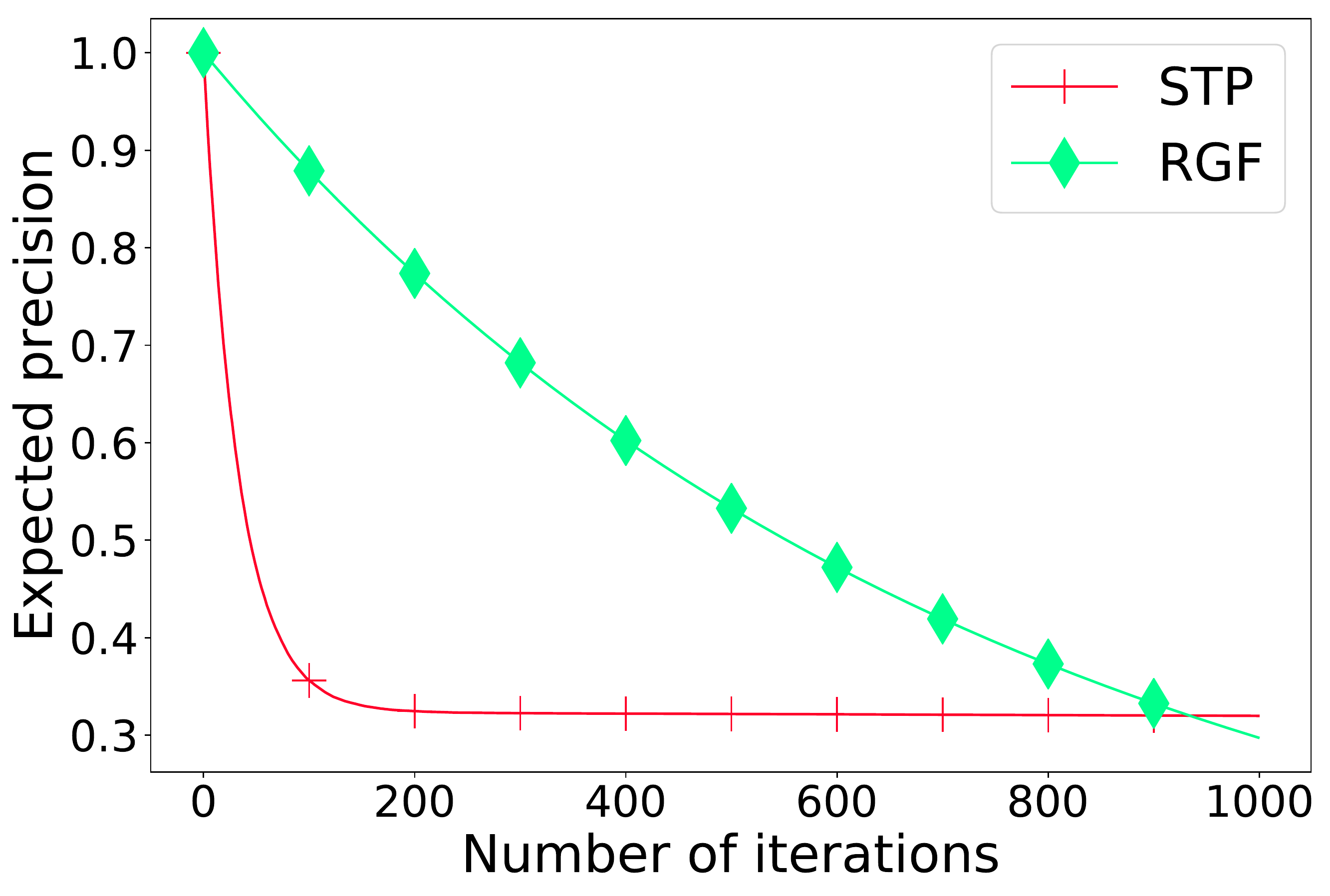}\label{subfig1:n50_rds_rgf}
	}
	\subfigure[$n = 75$]{
		\includegraphics[scale=0.145]{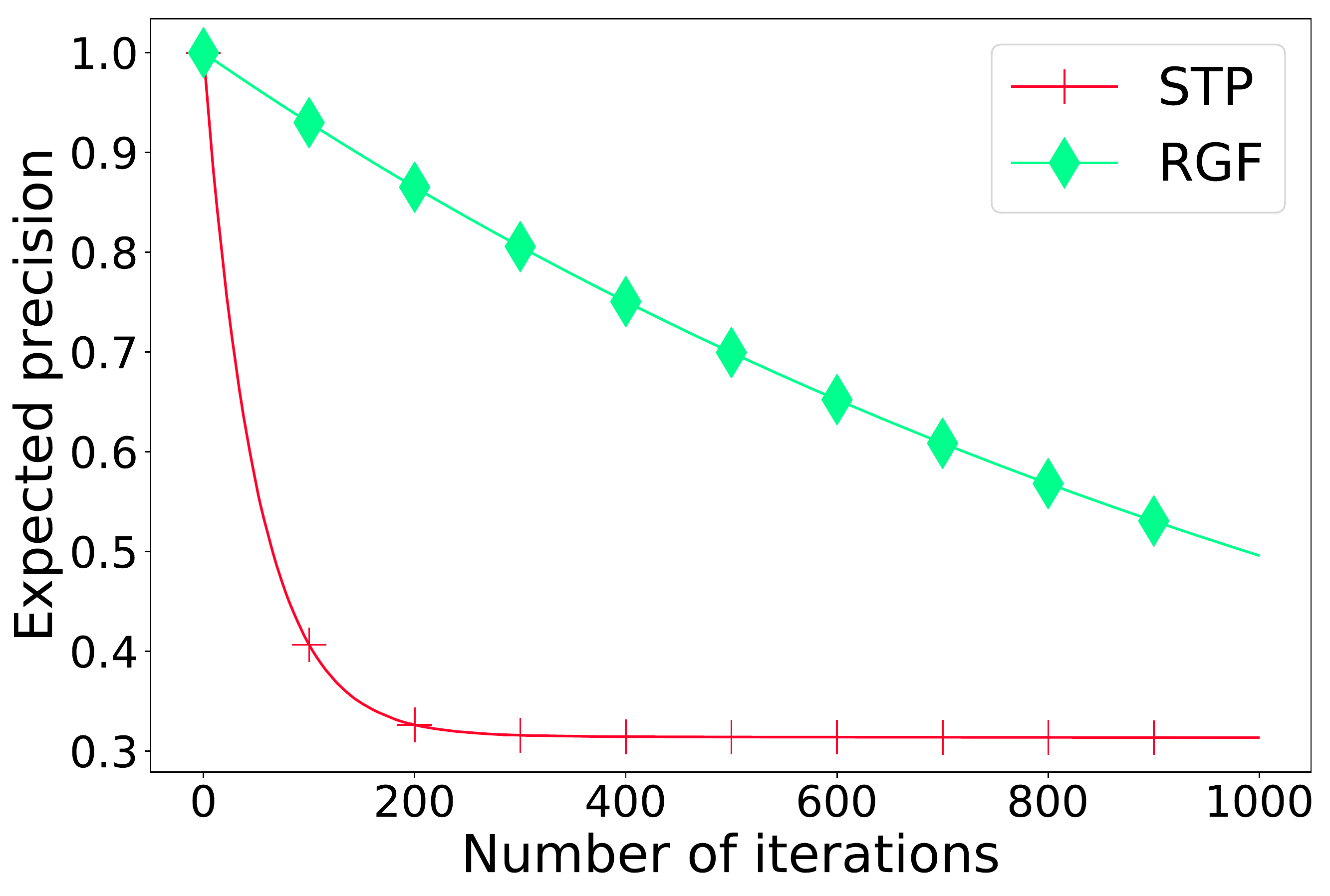}\label{subfig1:n75_rds_rgf}
	}
	\subfigure[$n = 100$]{
		\includegraphics[scale=0.145]{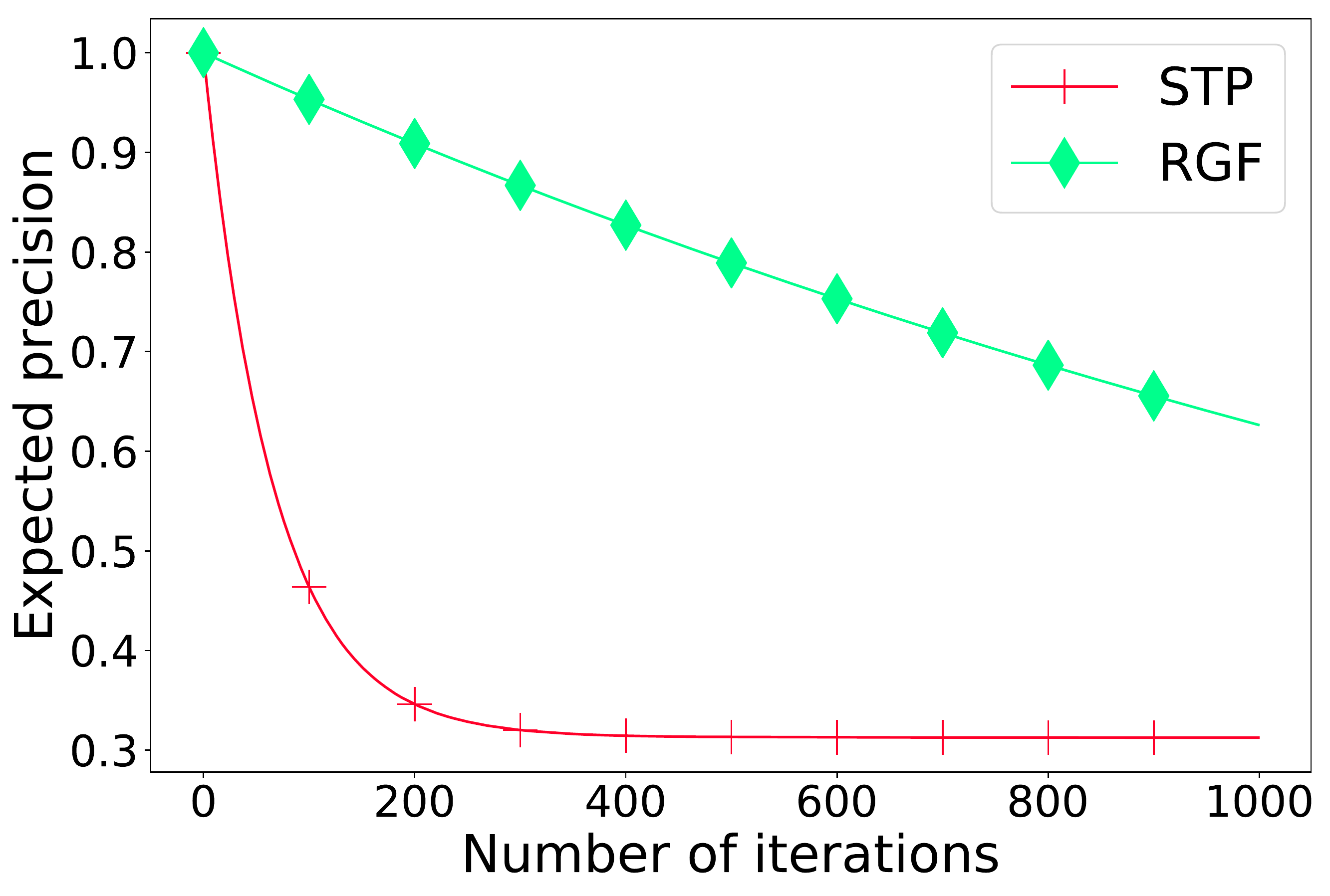}\label{subfig1:n100_rds_rgf}
	}
	\subfigure[$n = 125$]{
		\includegraphics[scale=0.145]{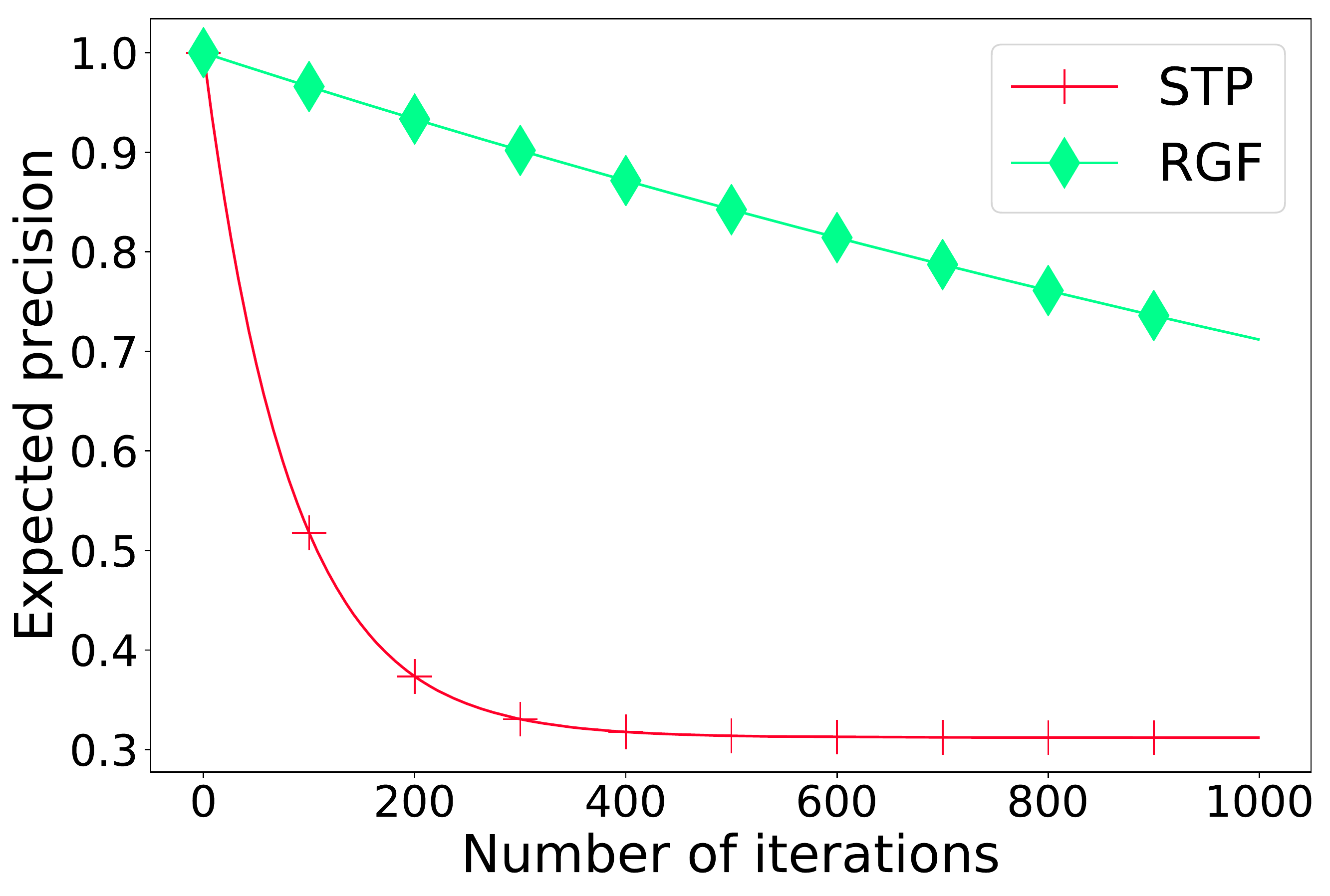}\label{subfig1:n125_rds_rgf}
	}
	\subfigure[$n = 150$]{
		\includegraphics[scale=0.145]{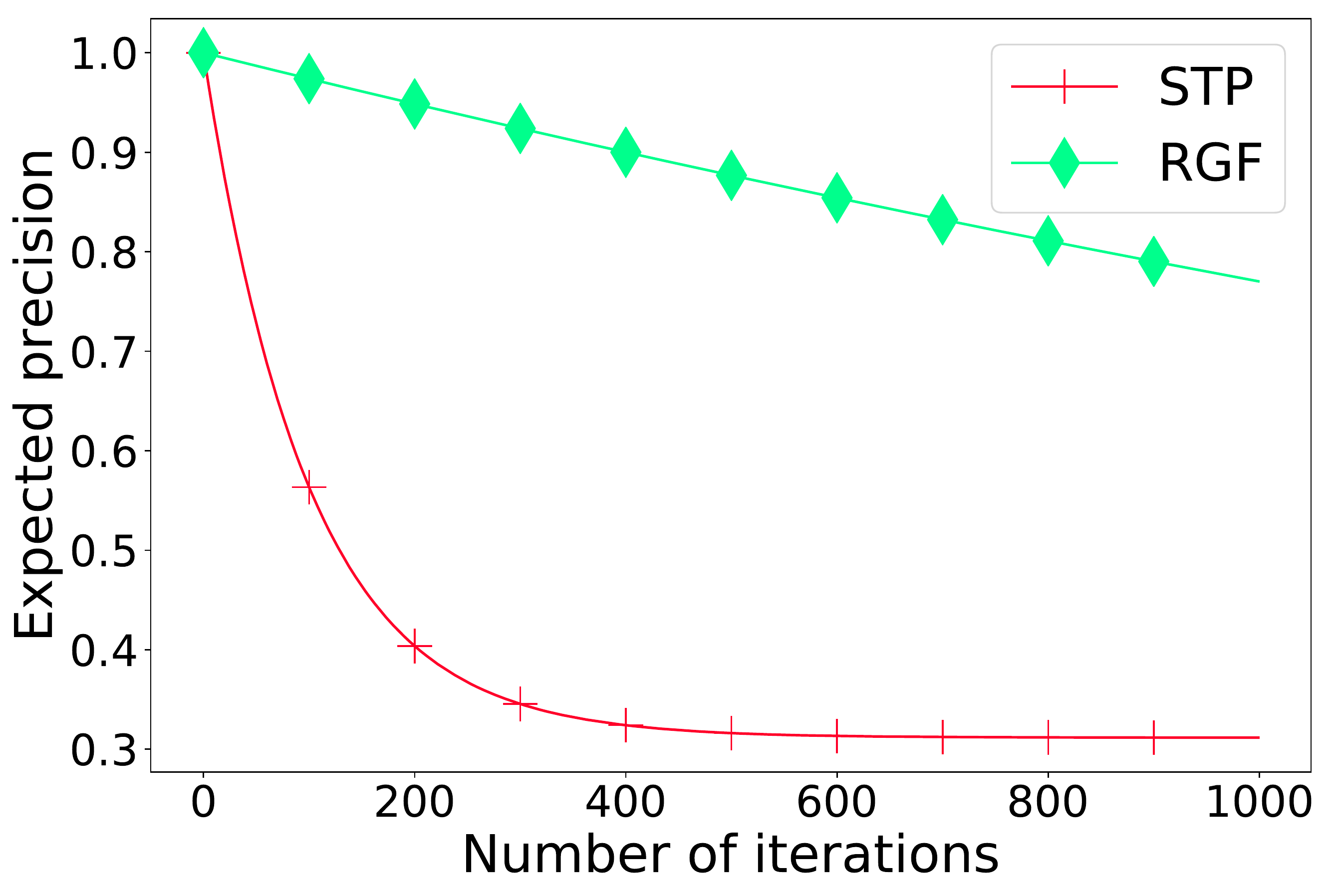}\label{subfig1:n150_rds_rgf}
	}
	\subfigure[$n = 175$]{
		\includegraphics[scale=0.145]{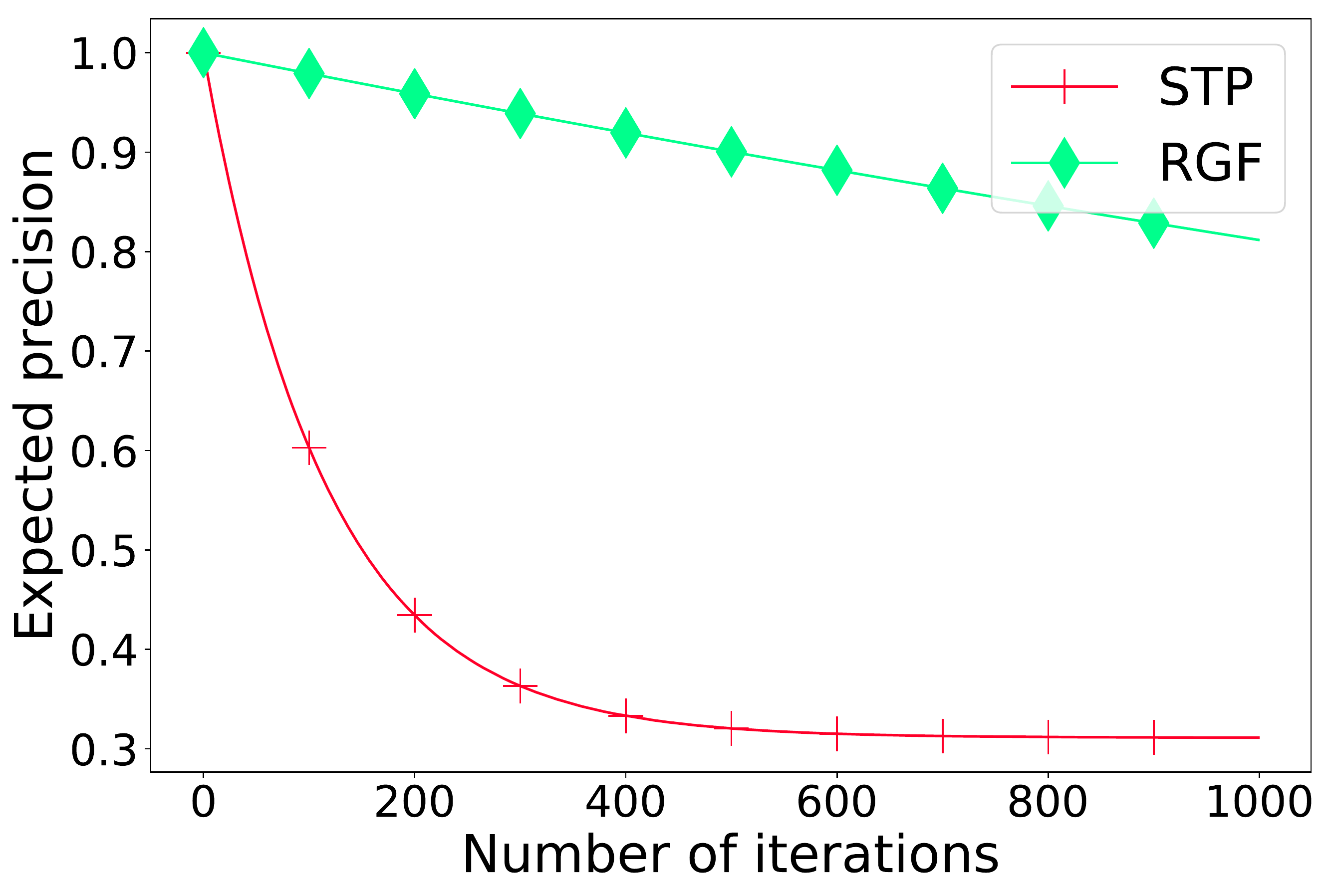}\label{subfig1:n175_rds_rgf}
	}
	\subfigure[$n = 200$]{
		\includegraphics[scale=0.145]{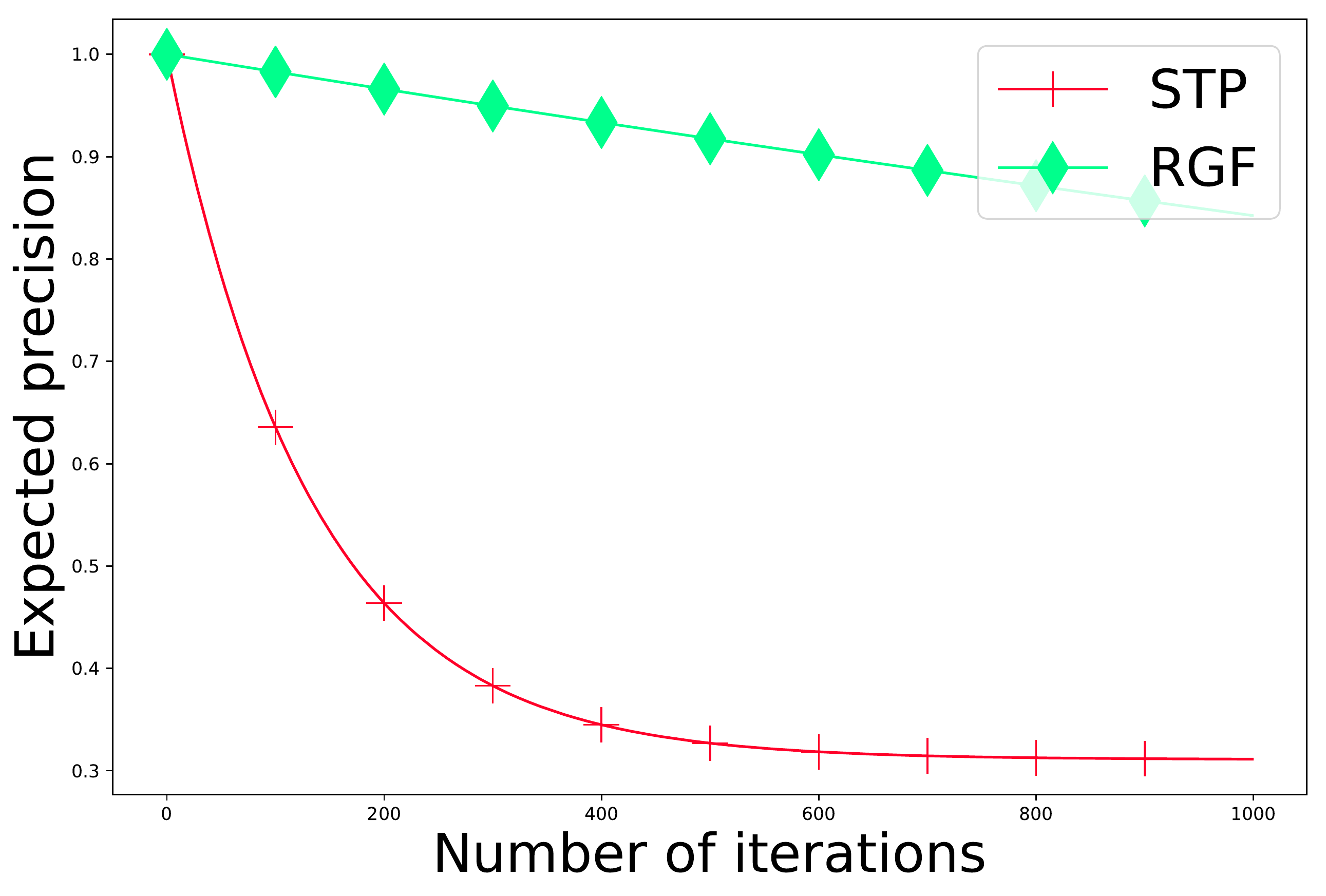}\label{subfig1:n200_rds_rgf}
	}
	\caption{Trajectories of \texttt{STP} and \texttt{RGF} for the different $n$.} \label{fig:_rds_rgf}
\end{figure}

We considered the following function
$$
f(x) = \frac{1}{2}x_1^2 + \frac{1}{2}\sum\limits_{i=1}^{n-1}(x_{i+1}-x_i)^2 + \frac{1}{2}x_n^2 - x_1
$$
and run \texttt{STP} and \texttt{RGF} for different $n$ (see Figure~\ref{fig:_rds_rgf}). We measure $\frac{f(x_k) - f_*}{f(x_0)-f_*}$ on the $y$-axis and call it ``Expected precision". One can notice that \texttt{STP} becomes more beneficial then \texttt{RGF} when $n$ is growing.

\section{Conclusions}
\label{sec:conc}

In this paper, we have proposed a very simple randomized algorithm --- Stochastic Three Points (\texttt{STP}) method ---  for derivative free optimization (DFO). At each iteration, the proposed method try to decrease the objective function 
along  a random direction sampled from a certain fixed probability law. 
Under mild assumption on this law, we have given the properties of this method for non-convex, convex and strongly convex problems. In fact, we have derived different practical rules for the stepsizes for which this method converges in expectation to a stationary point of the considered problem.

We have derived  the worst case complexity of \texttt{STP}. In fact, in the non-convex case, we have shown that \texttt{STP} needs  $O(n\e^{-2})$ function evaluations to find a point at which the gradient of the objective function is below $\e$, in expectation. 
In the convex case, the number of iterations to find a point at which the distance between the objective function and its optimal value in expectation is  $O(n \e^{-1})$. \texttt{STP} is shown to converge linearly for the strongly convex problems, i.e.\ the complexity is $O(n \log(1/\e))$. The complexity of \texttt{STP} depends linearly on the dimension of the considered problem, while this dependence is quadratic for deterministic direct search (\texttt{DDS}) methods.  We have also proposed a parallel version for \texttt{STP}.  

Our numerical experiments showed encouraging performance of the proposed \texttt{STP}  algorithm. 
A number of issues need further investigation, in particular the best choice of probability law for choosing the random directions. Extending our results to the non smooth problems and/or the constrained problems remains an interesting topic for the future research. 
It would be also interesting to confirm the potential of the proposed \texttt{STP} approach compared to the classical  approaches in DFO using extensive numerical tests.
\bibliographystyle{siam}
\bibliography{DFO-random}

\appendix
\section{Proof of Lemma~\ref{lem1}}
\label{app:A}
\begin{enumerate}
\item $\gamma_{\cal D} = \Exp \|s\|_2^2 = \tfrac{1}{A_n(1)}\int_{\|s\|_2^2=1} \|s\|_2^2 ds =  \tfrac{1}{A_n(1)}\int_{\|s\|_2^2=1} ds = 1$
where $A_n(1) = {\tfrac {2\pi ^{\tfrac {n}{2}}}{\Gamma ({\tfrac {n}{2}})}}$ is the air of the $n-1$-unit sphere and $\Gamma$ is the gamma function.

Let $\varepsilon_1 = g/\|g\|_2$ and $\varepsilon_2,\ldots,\varepsilon_n$ complete $\varepsilon_1$ to an orthonormal basis of $\R^n$ then
\begin{eqnarray*}
\Exp |\ve{g}{s}| &=&\tfrac{1}{A_n(1)} \int_{\|s\|_2^2=1}  |\ve{g}{s}| ds = \|g\|_2 \tfrac{1}{A_n(1)} \int_{\sum_{i=2}^n s_i^2=1 - s_1^2} |s_1|ds \\
&=& \|g\|_2 \tfrac{1}{A_n(1)} \int_{-1}^1 |s_1| \int_{\sum_{i=2}^n s_i^2=1 - s_1^2} ds_{2:n} ds_1\\
&=& \|g\|_2\tfrac{1}{A_n(1)}\int_{-1}^1 |s_1| A_{n-1}\left(1 - s_1^2\right) ds_1,
\end{eqnarray*}
where $A_{n-1}\left(1 - s_1^2\right) = \tfrac{2\pi^{(n-1)/2}\left(1 - s_1^2\right)^{n-2}}{\Gamma \left((n-1)/2\right)}$ is the volume of the $n-2$ sphere of radius $1 - s_1^2$, hence
\begin{eqnarray*}
\Exp |\ve{g}{s}| &=&  \|g\|_2  \tfrac{1}{A_n(1)}  \tfrac{2\pi^{(n-1)/2}}{\Gamma ((n-1)/2)} \int_{-1}^1 |s_1| \left(1 - s_1^2\right)^{n-2} ds_1 \\
&=&  \|g\|_2  \tfrac{1}{A_n(1)}  \tfrac{2\pi^{(n-1)/2}}{\Gamma ((n-1)/2)(n-1)}.
\end{eqnarray*}

If $n-1=2p$ then 
\begin{eqnarray*}
\Exp |\ve{g}{s}| &=&  \|g\|_2  \tfrac{2\pi^{p} \Gamma (p+1/2)}{ 2p\Gamma (p) 2 \pi^{p}\sqrt{\pi}} =  \|g\|_2   \tfrac{(2p)!}{2^{2p+1}(p!)^2} \quad \sim \quad \tfrac{\|g\|_2 }{ 2 \sqrt{\pi p}} ,
\end{eqnarray*}
since according to Stirling formula, $ p! \sim p^p e^{-p}\sqrt{2\pi p}$.  
If $n-1=2p+1$ then 
\begin{eqnarray*}
\Exp |\ve{g}{s}| &=&  \|g\|_2    \tfrac{2\pi^{p} \sqrt{\pi} \Gamma(p+1)} {2 \pi^{p+1} (2p+1)\Gamma (p+1/2)} =  
\|g\|_2   \tfrac{(p!)^2 2^{2p}}{(2p+1)!\pi} \sim \tfrac{ \sqrt{p}} {\sqrt{\pi}(2p+1)} \sim \tfrac{\|g\|_2 }{ 2 \sqrt{\pi p}} 
\end{eqnarray*}

In the both cases, $\Exp |\ve{g}{s}| \sim  \tfrac{\|g\|_2 }{ 2 \sqrt{\pi p}} \sim  \tfrac{\|g\|_2 }{  \sqrt{2\pi n}}.$
\item $\gamma_{\cal D} = \Exp \|s\|_2^2 = \tfrac{1}{n}\Exp \|x\|_2^2 =1$, where $x \sim N(0,I)$.

Note that $s \sim \tfrac{1}{\sqrt{n}} N(0,I)$ implies $ \ve{g}{s} \sim \tfrac{1}{\sqrt{n}} N(0,\|g\|_2^2)$, hence
\begin{eqnarray*}
\Exp |\ve{g}{s}| &=& \tfrac{1}{\|g\|_2  \sqrt{2n\pi}}\int_{-\infty}^{+\infty} |x| e^{-\tfrac{x^2}{2\|g\|_2^2}} dx 
= \tfrac{\sqrt{2}}{  \sqrt{n\pi}}\|g\|_2.
\end{eqnarray*}

\item $\gamma_{\cal D} =  \sum_{i=1}^n \|e_i\|_2^2 P(s=e_i) =1$ and $\Exp |\ve{g}{s}|  = \tfrac{1}{n}\sum_{i=1}^n |g_i| = \tfrac{1}{n}\|g\|_1.$

\item  $\gamma_{\cal D} =  \sum_{i=1}^n \|e_i\|_2^2 P(s=e_i) =1$ and
$\Exp |\ve{g}{s}|  = \sum_{i=1}^n |g_i| P(s=e_i) = \sum_{i=1}^n p_i |g_i|.$

\item  $\gamma_{\cal D} = \sum_{i=1}^n \|d_i\|_2^2 P(s = d_i) = \sum_{i=1}^n p_i=1$ and 
$
\Exp |\ve{g}{s}| = \sum_{i=1}^n p_i |g_i d_i | = \|g\|_{\cal D}.
$
\end{enumerate}

\section{Proof that our approach covers some first order methods}
\label{app:B}

\begin{itemize}
\item Normalized Gradient Descent (NGD) method: 

At iteration $k$, $s \sim {\cal D}_k$ means that 
$s = \tfrac{g_k}{\|g_k\|_2}$ with probability 1.
\begin{eqnarray*}
	\gamma_{{\cal D}_k} = \Exp_{s\sim {\cal D}_k} \|s\|_2^2 = 1,\\
	\Exp_{s\sim {\cal D}_k}\; |\ve{g_k}{s}|  = \|g_k\|_2.
\end{eqnarray*}


\item Signed Gradient Descent (SignGD) method:

At iteration $k$, $s \sim {\cal D}_k$ means that $s = sign \left( g_k \right)$ with probability $1$, where the $sign$ operation is element wise sign.
\begin{eqnarray*}
\gamma_{{\cal D}_k} =  \Exp_{s\sim {\cal D}_k} \|s\|_2^2 = \Exp_{s\sim {\cal D}_k} \| sign \left( {g}_k \right) \|_2^2 \le \sum_{i=1}^n 1 = n,\\
\Exp_{s\sim {\cal D}_k}\; |\ve{g_k}{s}|  = \Exp_{s\sim {\cal D}_k}\; |\ve{g_k}{sign \left( {g}_k \right)}| = \|g_k\|_1.
\end{eqnarray*}

%

\item  Normalized Randomized Coordinate Descent (NRCD) method (equivalently this method can be called Randomized Signed Gradient Descent): 

At iteration $k$, $s \sim {\cal D}_k$ means that 
$s = \tfrac{g_k^{i}}{ |g_k^{i} |}e_i$ with probability $\tfrac{1}{n}$, where $g_k^{i}$ is the $i-th$ component of $g_k$. 
\begin{eqnarray*}
\gamma_{{\cal D}_k} =  \Exp_{s\sim {\cal D}_k} \|s\|_2^2 = \tfrac{1}{n} \sum_{i=1}^n 1  = 1\\
\Exp_{s\sim {\cal D}_k}\; |\ve{g_k}{s}|  = \Exp_{i \sim {U[1,\ldots,n]}}\; \left|\ve{g_k}{\tfrac{g_k^{i}}{ |g_k^{i} |}e_i}\right| = \tfrac{1}{n} \sum_{i=1}^n  |g_k^{i} |   = \tfrac{1}{n} \|g_k\|_1.
\end{eqnarray*}

\item Normalized Stochastic Gradient Descent (NSGD) method:

At iteration $k$, $s \sim {\cal D}_k$ means that 
$s = \hat{g}_k $ where $\hat{g}_k$ is 
 the  stochastic gradient satisfying $\Exp\left[ \hat{g}_k \right] = \tfrac{g_k}{\|g_k\|_2}$, and $\Exp\left[ \|\hat{g}_k\|_2^2  \right]  \le \sigma < \infty $.
\begin{eqnarray*}
\Exp_{s\sim {\cal D}_k}\; |\ve{g_k}{s}|  &=& \Exp_{s\sim {\cal D}_k}\; |\ve{g_k}{ \hat{g}_k }|  \ge  \Exp_{s\sim {\cal D}_k}\; \ve{g_k}{ \hat{g}_k} = \|g_k\|_2.
\end{eqnarray*}
\end{itemize}
\end{document}